\input amstex
\documentstyle{amsppt}
\magnification1200
\tolerance=10000
\def\n#1{\Bbb #1}
\def\p{\Bbb C_{\infty}}

\def\Exp{\hbox{Exp}}

\def\End{\hbox{End}}

\def\Lie{\hbox{Lie}}

\def\ord{\hbox{ord }}
\def\max{\hbox{max }}
\def\Lie{\hbox{Lie }}

\def\e11{E_{11}}

\def\ga{\goth A}

\def\ve{\varepsilon}
\def\de{\delta}
\def\De{\Delta}
\def\ga{\gamma}

\def\be{\beta}
\def\th{\theta}
\def\al{\alpha}
\def\om{\omega}
\def\la{\lambda}
\def\vf{\varphi}

\def\g{\goth }

\topmatter
\title
Lattice map for Anderson t-motives: first approach
\endtitle
\author
A. Grishkov, D. Logachev
\endauthor
\NoRunningHeads
\address
First author: Departamento de Matem\'atica e estatistica
Universidade de S\~ao Paulo. Rua de Mat\~ao 1010, CEP 05508-090, S\~ao Paulo, Brasil, and Omsk State University n.a. F.M.Dostoevskii. Pr. Mira 55-A, Omsk 644077, Russia.
Second author: Departamento de Matem\'atica, Universidade Federal de Amasonas, Manaus, Brasil
\endaddress
\thanks Thanks: The authors are grateful to FAPESP, S\~ao Paulo, Brazil for a financial support (process No. 2013/10596-8). The first author is grateful to SNPq, Brazil, to RFBR, Russia, grant 16-01-00577a (Secs. 1-4), and to Russian Science Foundation, project 16-11-10002 (Secs. 5-8) for a financial support. The second author is grateful to Max-Planck-Institut f\"ur Mathematik, Bonn; Hausdorff Institut f\"ur Mathematik, Bonn; to Sinnou David, R\'egis de la Br\`eteche, Lo\"{\i}c Merel, Institut de Math\'ematiques de Jussieu, UPMC, Paris, for invitations in academic years 2012--13, 2013--14. Also, the authors are grateful to an anonimous reviewer for careful reading of the manuscript and providing a list of mistakes.
\endthanks
\abstract There exists a lattice map from the set of pure uniformizable Anderson t-motives to the set of lattices. It is not known what is the image and the fibers of this map. We prove a local result that sheds the first light to this problem and suggests that maybe this map is close to 1 -- 1. Namely, let $M(0)$ be a t-motive of dimension $n$ and rank $r=2n$ \ --- \ the $n$-th power of the Carlitz module of rank 2, and let $M$ be a t-motive which is in some sense "close" to $M(0)$. We consider the lattice map $M \mapsto L(M)$, where $L(M)$ is a lattice in $\p^n$. We show that the lattice map is an isomorphism in a "neighborhood" of $M(0)$. Namely, we compare the action of monodromy groups: (a) from the set of equations defining t-motives to the set of t-motives themselves, and (b) from the set of Siegel matrices to the set of lattices. The result of the present paper gives that the size of a neighborhood, where we have an isomorphism, depends on an element of the monodromy group. We do not know whether there exists a universal neighborhood. Method of the proof: explicit solution of an equation describing an isomorphism between two t-motives by a method of successive approximations using a version of the Hensel lemma.
\endabstract
\endtopmatter
\document
{\bf 0. Introduction.}
\nopagebreak
\medskip
t-motives ([G], 5.4.2, 5.4.18, 5.4.16) are the functional field analogs of abelian varieties (more exactly, of abelian varieties with multiplication by an imaginary quadratic field, see for example [L1]). For the number field case we have a classical theorem (here and below we consider lattices up to a linear transformation of the ambient space):
\medskip
{\bf Theorem 0.1.} Abelian varieties of dimension $g$ over $\n C$ are in 1 -- 1 correspondence with lattices satisfying Riemann condition, of dimension $2g$ in $\n C^g$.
\medskip
Our knowledge on the functional field analog of this theorem is very poor, and the purpose of the present paper is to get a result towards this analog. Let $\p$ be the analog of $\n C$ in characteristic $p$ (it is a complete algebraically closed field). Throughout all the paper we consider for simplicity only t-motives $M$ over the affine line $A^1$ such that their nilpotent operators $N$ (see (1.3.1) below for its definition) are equal to 0. Let $r$, $n$ be respectively the rank and dimension of $M$. For some $M$ it is possible associate to $M$ a lattice $L(M)$ of rank $r$ in $n$-dimensional space $\p^n$ (see (1.4) for a definition of a lattice). These $M$ are called uniformizable ([A], Section 2). $M \mapsto L(M)$ is a contravariant functor.
\medskip
For $n=1$ the situation is completely analogous to the Theorem 0.1:
\medskip
{\bf Theorem 0.2} ([Dr]). All t-motives of dimension 1 ( = Drinfeld modules) are uniformizable. There is a 1 -- 1 correspondence between Drinfeld modules of rank $r$ over $\p$ and lattices of rank $r$ in $\p$.
\medskip
There exists a notion of purity of $M$ (see [G], 5.5.2 for the definition; all Drinfeld modules are pure). For $n=r-1$ the duality theory gives us an immediate corollary of Theorem 0.2:
\medskip
{\bf Corollary 0.3} ([L], Corollary 8.4). All pure t-motives of rank $r$ and dimension $r-1$ over $\p$ are uniformizable. There is a 1 -- 1 correspondence between their set, and the set of lattices of rank $r$ in $\p^{r-1}$ having dual.
\medskip
Not all such lattices have dual, but almost all, i.e. even in this simple case the correspondence is not strictly 1 -- 1, but only an "almost 1 -- 1".
\medskip
For arbitrary $n$, $r$ we have
\medskip
{\bf Theorem 0.4} ([A]). If $M$ is uniformizable then its lattice $L(M)$  is well-defined. Not all $M$ are uniformizable.
\medskip
We know neither the image of the lattice map $M \mapsto L(M)$ nor its fibre. Taking into consideration
\medskip
{\bf Theorem 0.5} ([H], Theorem 3.2). The dimension of the moduli space of pure t-motives of rank $r$ and dimension $n$ is equal to $n(r-n)$.
\medskip
and the obvious fact that the moduli space of lattices of rank $r$ in $\p^n$ has the same dimension $n(r-n)$ we can state
\medskip
{\bf Conjecture 0.6.} Let us consider the lattice map $M \mapsto L(M)$ from the set of pure uniformizable t-motives to the set of lattices. Its image is open, and its fibre at a generic point is discrete.
\medskip
{\bf Remark 0.6.1.} Both the set of pure uniformizable t-motives and the set of lattices are quotient sets of some sets of matrices (the matrices entering in equations defining a t-motive, and Siegel matrices of lattices). We consider in Conjecture 0.6 the quotient topologies. Apparently it is not known much on these topologies, for example, there is no proof of existence of fundamental domains. Proposition 1.7.2 shows that $\om I_n$ (see below for $\om$) is an isolated point of the orbit of $GL_{2n}(\n F_q[\theta])$, while Proposition 1.7.1 shows that the action of $GL_{2n}(\n F_q[\theta])$ on the set of Siegel matrices is not as good as the corresponding action in the number field case.
\medskip
{\bf Remark 0.7.} Preliminary results of [L2] suggest that the condition of purity in 0.6 is essential: the dimension of the fibre conjecturally can be $>0$ for the non-pure case.
\medskip
Theorem 0.9 --- the main result of the  present paper --- is a local result supporting 0.6. Let $q$ be a power of $p$, $\theta$ a transcendent element and $\n F_q[\theta]\subset \p$ the function field case analog of $\n Z$. Further, let $A\in M_n(\p)$. We consider a t-motive $M(A)$ (see below for the details) defined as follows: $M(A)$ is a module over the Anderson ring $\p[T,\tau]$ (see (1.1)) which is free of dimension $n$ over the ring $\p\{\tau\}$, and for its $\p\{\tau\}$-basis $e_*=(e_1,\dots,e_n)^t$ (written as a column) the action of multiplication by $T$ is given by the following explicit equation:
$$Te_*=\theta e_*+A\tau e_* +\tau^2e_*\eqno{(0.8)}$$
The rank of $M(A)$ is $2n$. The t-motive $M(0)$ is the initial object of the present paper, we shall work in its neighborhood. We have $M(0)=\g C_2^{\oplus n}$ where $\g C_2$ is the Carlitz module over $\n F_{q^2}$.

There exists a neighborhood of 0 such that if $A$ belongs to it then $M(A)$ is uniformizable. Very roughly speaking, we give some evidence that Conjecture 0.6 holds, and moreover that the lattice map is a 1 -- 1 correspondence in a "system of neighborhoods" of $M(0)$.
\medskip
Let us formulate the theorem and outline the proof.  A definition of a Siegel matrix of a lattice of dimension $n$ and rank $2n$ for the function field case is the same as in the number field case (see 1.5). It is a $n\times n$ matrix with entries in $\p$. A Siegel matrix $S$ defines its lattice $\goth L(S)$ uniquely, while a lattice has many Siegel matrices: like in the number field case, the group $GL_{2n}(\n F_q[\theta])$ (almost, see 1.6.4) acts on the set of Siegel matrices, and two Siegel matrices $S_1$, $S_2$ define the same lattice if and only if $\exists \gamma\in GL_{2n}(\n F_q[\theta])$ such that $\gamma(S_1)=S_2$.
\medskip
For simplicity, we shall consider only the case of odd $q$, except Proposition 1.7.1 (the case of even $q$ requires minor modifications). Let us fix throughout the whole paper an element $\omega\in \n F_{q^2}-\n F_{q}$ such that $\om^2\in \n F_{q}$. A Siegel matrix of $L(M(0))$ is $\omega I_n$. We consider the group $G_0\subset GL_{2n}(\n F_q[\theta])$ --- the stabilizer of $\omega I_n$, i.e. the monodromy group of the map (Siegel matrices) $\to$ (lattices) at the lattice $L(M(0))$. This group is in some sense the "biggest" among monodromy groups of other elements of the set of Siegel matrices (this explains why we consider neighborhood of $M(0)$ but not of other t-motive). It is isomorphic to $GL_n(\n F_{q^2}[\theta])$.

For a sufficiently small $A$ it is possible to choose the distinguished representative $\Cal S(A)$ in the set of all Siegel matrices of $L(M(A))$ (i.e. $\goth L(\Cal S(A))=L(M(A))$ ) satisfying the condition that $\Cal S(A)$ is close to $\omega I_n$, see 4 lines above (2.10) for the definition of $\Cal S(A)$. Entries of $\Cal S(A)$ are power series of the entries of $A$, and we can consider $\Cal S$ as a map from a neighborhood of 0 in $M_n(\p)$ to a neighborhood of $\omega I_n$ in $M_n(\p)$. Since the action of $\gamma\in GL_{2n}(\n F_q[\theta])$ is continuous, we have: if $S\in M_n(\p)$ is close to $\omega I_n$ and $\gamma\in G_0$ then $\gamma(S)$ is close to $\omega I_n$.
\medskip
{\bf Theorem 0.9.} (1). Surjectivity. $\Cal S$ is 1 -- 1 from a neighborhood of 0 to a neighborhood of $\omega I_n\in M_n(\p)$ --- the set of Siegel matrices. Particularly, the lattice map $M\mapsto L(M)$ is surjective near the lattice $L(M(0))$.
\medskip
(2). Injectivity.  For any $m\ge 0$ there exists $U_m$ --- a neighborhood of 0 in $M_n(\p)$ having the following property. Let $A_1$, $A_2\in U_m$ and $M(A_1)$, $M(A_2)$ the corresponding t-motives. If $\exists \gamma\in G_0$ such that the entries of $\gamma$, $\ga^{-1}$ have degrees (as polynomials in $\theta$) \  $\le m$ and $\gamma(\Cal S(A_1))=\Cal S(A_2)$ then $M(A_1)$ is isomorphic to $M(A_2)$.
\medskip
{\bf Remark 0.9.A.} We have $$L(M(A_1))=L(M(A_2))\Longrightarrow \exists \ \ga \in GL_{2n}(\n F_q[\theta]) \hbox{ such that }\gamma(\Cal S(A_1))=\Cal S(A_2)\eqno{(0.9.A.1)}$$ hence (0.9.2) practically means that if $\ga$ of (0.9.A.1) satisfies conditions of (0.9.2) then $$\{L(M(A_1))=L(M(A_2))\}\Longrightarrow \{M(A_1) \hbox{ is isomorphic to }M(A_2)\}$$ i.e. (0.9.2) really indicates to injectivity of the map $M\mapsto L(M)$.
\medskip
{\bf Remark 0.9.B.} Let $A_1$, $A_2\in U_m$, where $m>>0$ is fixed, and $\gamma\in GL_{2n}(\n F_q[\theta])$ such that $\gamma(\Cal S(A_1))=\Cal S(A_2)$. In general, $\ga\not\in G_0$, see Proposition 1.7.1. There is a problem to find an analog of Theorem 0.9 for this case, see 0.10.
\medskip
{\bf Method of the proof.} To prove (0.9.1) we show that $\Cal S(A)$ is a power series of $A$. A version of the Hensel lemma (Lemma 2.29) shows that this series is 1 -- 1 in a neighborhood of 0. This is made in Section 2.

To prove (0.9.2) we consider equalities in the matrix ring $M_n(\p\{\tau\})$. We denote $$T_A:=\theta +A\tau+\tau^2\in M_n(\p\{\tau\})\eqno{(0.9.3)}$$ --- a linear transformation of $\p\{\tau\}^n$. Particularly, $T_0:=\theta+\tau^2$. To prove that $M(A_1)$ is isomorphic to $M(A_2)$ it is sufficient to find $\hat \goth B_2\in GL_n(\p\{\tau\})\subset M_n(\p\{\tau\})$ such that $$T_{A_1}\hat \goth B_2=\hat \goth B_2T_{A_2}\eqno{(0.9.4)}$$
(equality in $M_n(\p\{\tau\})$ ). Since $A_1$, $A_2\approx 0$ we shall find $\hat \goth B_2$ such that $\hat \goth B_2\approx \goth B_2$ where $\goth B_2$ has the property that $$T_{0}\goth B_2=\goth B_2T_{0}\eqno{(0.9.5)}$$
The set of $\goth B_2$ satisfying (0.9.5) is denoted by $G_2$. There is a canonical isomorphism $\al\circ\be:G_0\to G_2$. We choose $\goth B_2=\al\circ\be(\ga)$ and $\hat \goth B_2=\goth B_2+Y$ where $Y\approx0$. Further, in (0.9.4) we consider $A_1$ as a parameter, $Y$ and $A_2$ as unknowns. (0.9.4) becomes a system of matricial equations. It is necessary to emphasize that at the first glance it seems that (0.9.4) has no solutions satisfying $\hat \goth B_2\not\in GL_n(\p)$. Nevertheless, such solutions really exist, and moreover it turns out that the system (0.9.4) can be solved by a method of successive approximations, using a version of the Hensel lemma (Lemma 2.29), if $A_1$ is sufficiently small. This is proved in Proposition 4.2. More exactly, we write (0.9.4) in the form $T_{A_1}\hat \goth B_2=\hat \goth B_2T_X$ where $X$ is an unknown matrix. We show that $X$ exists; now we must prove that $X$ obtained as a solution to (0.9.4) satisfies $X=A_2$. This is made in Lemma 4.23 and in 4.24.
\medskip
{\bf Remark 0.9.6.} $G_0$ acts on the set of Siegel matrices in a  neighborhood of $\omega I_n\in M_n(\p)$, while (0.9.4) shows that - strictly speaking - $G_2$ does not act on the set of $T_A$, where $A\approx 0$, but only a "modification" $\hat \goth B_2$ of $\goth B_2\in G_2$, and this modification $\hat \goth B_2$ depends on $A$. It would be interesting to axiomatize this phenomenon.
\medskip
{\bf 0.10. Further research.} The final purpose of the present research is to prove or to disprove Conjecture 0.6, to find the image and the fibres of the lattice map. According Proposition 1.7.1, it can happen that $\ga$ from 0.9, (2) does not belong to $G_0$. Namely, let $S_1$, $S_2$ be near to $\om I_n$ matrices and $\ga \in GL_{2n}(\n F_q[\th])$, $\ga\not\in G_0$ such that $S_2=\ga(S_1)$. We consider the corresponding near to 0 matrices $\Cal S^{-1}(S_1)$, $\Cal S^{-1}(S_2)$ and the corresponding t-motives $M_1:=M(\Cal S^{-1}(S_1))$, $M_2:=M(\Cal S^{-1}(S_2))$. Because of $S_2=\ga(S_1)$, the lattices of $M_1$, $M_2$ are isomorphic.
\medskip
{\bf Problem 0.10.1.} Are $M_1$, $M_2$ isomorphic?
\medskip
This case is not covered by Theorem 0.9. If the answer is yes at least for a one non-trivial case, this gives much more evidence in favor of the conjecture that the lattice map is 1 -- 1 near $M(0)$.
\medskip
{\bf Idea of a solution.} We need to find $\hat \g B_2$ satisfying (0.9.4). We cannot choose $\g B_2$ satisfying (0.9.5) as a first approximation to $\hat \g B_2$, because $\ga\not\in G_0$. Alternatively, we can consider $M_1$, $M_2$ as $\p[T]$-modules and to find an analog of (0.9.4) for $\p[T]$-modules:
$$Q_1\hat \g B_T=\hat \g B^{(1)}_TQ_2\eqno{(0.10.2)}$$ where $Q_1$, resp. $Q_2$ are matrices of multiplication by $\tau$ of $M_1$, $M_2$ treated as $\p[T]$-modules, $\hat \g B_T\in GL_4(\p[T])$ is a matrix of a $\p[T]$-isomorphism between them, and for $P=\sum a_iT^i\in \p[T]$ we denote $P^{(1)}:=\sum a_i^qT^i$. Apparently it is difficult to find a relation between a $\g B_T$ --- a first approximation to $\hat \g B_T$ --- and $\ga$, if $\ga\not\in G_0$.
\medskip
So, we can try to find first approximations $\g B_2$ or $\g B_T$ by a computer search. We can assume that their entries are polynomials in $\tau$ or $T$ of small degree (probably of degree 1). The next step of solution: to show that (one of) these first approximations can be deformed to the exact solution of (0.9.4) or (0.10.2), using the methods similar to the ones of the present paper.
\medskip

For the example of Proposition 1.7.1 and for its small deformations the situation is similar to the one of the present paper: we can find $\g B_2$. Really, let $\ga$ and $S$ be from 1.7.1. Let $S_1$, $S_2=\ga(S_1)$ belong to a neighborhood of $S$. The lattice $\g L(S)$ is a direct sum of two (isogenous, with complete multiplication) lattices of dimension 1 and rank 2. Hence, the corresponding t-motive $M$ is a direct sum of two Drinfeld modules of rank 2, and the lattice functor gives us an isomorphism $i: \End(M)\to\End(\g L(S))$. Hence, $\exists \ \vf\in \End(M)$ such that $i(\vf)=\ga$ of (1.7.1.1), and the method of the present paper can be used for this case: the map $\vf$ can be taken as the first approximation to the isomorphism between $M(\Cal S^{-1}(S_1))$, $M(\Cal S^{-1}(S_2))$. Clearly for this case $M(\Cal S^{-1}(S))$ plays the same role as $M(0)$ for the present paper. 
\medskip
Does exist a less trivial counterexample of Proposition 1.7.1? This is a subject of further research. 
\medskip
{\bf 0.10.3.} As the next step, we should try to find a universal $U_0$ --- a neighborhood of 0 in $M_n(\p)$ --- such that 0.9 holds for this $U_0$ for all $\ga$ (the present proof of (0.9) gives a rapidly decreasing sequence of $U_m$). We have an obvious
\medskip
{\bf Corollary 0.11.} Let $U_1$, $A_1=(a_{1ij})$, $A_2=(a_{2ij})$, $m$, $\ga$ be as in (0.9.2), and $W_0$ a number such that if ord $a_{ij}\ge W_0$ (see beginning of Section 1 for the definition of the function ord) then $A\in U_1$. If ord $a_{1ij}$, ord $a_{2ij}\ge W_0+2m$ and $\ga$ is a product $\ga=\ga_1\cdot \dots \cdot \ga_m$ where $\ga_i\in G_0$ are such that the entries of $\gamma_i$, $\ga_i^{-1}$ have degrees (as polynomials in $\theta$) \  $\le 1$, then (0.9(2)) holds for these $A_1$, $A_2$, $\ga$.
\medskip
{\bf Deduction  from 0.9.} We define $A'_i$ by the formula $\gamma_i\cdot\gamma_{i-1}\cdot\dots\cdot\gamma_1(\Cal S(A_1))=\Cal S(A'_i)$, hence $A'_0=A_1$, $A'_m=A_2$. We have ord $(A'_i)_{jk}\ge W_0+2m-2i$ (see 4.2(4)), and (0.9.(2)) implies that $M(A'_i)$ is isomorphic to $M(A'_{i+1})$. $\square$
\medskip
A significant part (but clearly not all) of $\ga\in G_0$ can be represented as a product $\ga=\ga_1\cdot \dots \cdot \ga_m$. Hence, this corollary gives a hope to find a universal $U_0$.
\medskip
{\bf 1. Definitions. }

\nopagebreak
\medskip
Let $q$ be a power of a prime $p$. The field $\n F_q(\theta)$ is the functional field analog of $\n Q$. It has a valuation function ord: $\n F_q(\theta)^*\to \n Z$ defined by $\ord(f)=$ minus degree of $f$, where $f\in \n F_q(\theta)^*$ is a rational function. The ring $\n F_q[\theta]$ is the functional field analog of $\n Z$. The completion of $\n F_q(\theta)$ with respect to the topology defined by the valuation ord is the field of the Laurent series $\n F_q((1/\theta))$ --- the functional field analog of $\n R$. By definition, $\p$ is the completion of its algebraic closure. The valuation function ord can be prolonged uniquely to $\p$. For any matrix $A=(a_{ij})$ with entries $a_{ij} \in \p$ we define $\ord(A)=$ min ord $a_{ij}$, this gives a topology on the set of $A$.
\medskip
Let $\p[T,\tau]$ be the Anderson ring, i.e. the ring of non-commutative polynomials satisfying the following relations (here
$a \in \p$):
$$Ta=aT, \ T\tau = \tau T, \ \tau a = a^q \tau \eqno{(1.1)}$$
\medskip
{\bf Definition 1.2.} ([G], 5.4.2, 5.4.18, 5.4.16). A t-motive\footnotemark \footnotetext{Terminology of Anderson; Goss calls these objects abelian t-motives.} $M$ is a left
$\p[T, \tau]$-module which is free and finitely generated as both $\p[T]$-,
$\p\{\tau\}$-module and such that
$$ \exists \goth m= \goth m(M) \ \hbox{such that}\ (T-\theta)^ \goth m M/\tau M=0\eqno{(1.2.1)}$$
\medskip
The dimension of $M$ over $\p\{\tau\}$ (resp. $\p[T]$) is denoted by $n$ (resp.
$r$), these numbers are called the dimension and rank of $M$.
\medskip
We shall need the explicit matrix description of t-motives. Let
$e_*=(e_1, ..., e_n)^t$ be the vector column of elements of a basis
of $M$ over $\p\{\tau\}$. There exists a matrix $\goth A\in M_n(\p\{\tau\})$ such that

$$T e_* = \goth A e_*, \ \ \goth A = \sum_{i=0}^l \goth A_i \tau^i \hbox{ where } \goth A_i
\in M_n(\p)\eqno{(1.3)}$$
Condition (1.2.1) is equivalent to the condition
$$\goth A_0=\theta I_n + N\eqno{(1.3.1)}$$
where $N$ is a
nilpotent matrix, and the
condition \{$ \goth m(M)$ can be taken to 1\} is equivalent to the condition $N=0$.
\medskip
We fix $n$, and we shall consider only those $M$ whose equation (1.3) has the form (0.8), or, the same, $\goth A =T_A$ from (0.9.3), i.e. $l=2$, $N=0$, $\goth A_1=A$, $\goth A_2=I_n$. They have dimension $n$, rank $2n$, they are all pure, and there exists a neighborhood of 0 in $M_n(\p)$ such that if $A$ belongs to it then $M$ is uniformizable.
\medskip
{\bf Definition 1.4.} Let $V$ be the space $\p^n$. A free $r$-dimensional
$\n F_q[\theta]$-submodule
$L$ of $V$ is called a lattice if
\medskip
(a) $L$ generates $V$ as a $\p$-module and
\medskip
(b) The $\n F_q((1/\theta))$-linear span of $L$ has dimension $r$ over
$\n F_q((1/\theta))$.
\medskip
Two lattices $L_1,$ $L_2\subset V$ are called isomorphic if there exists a $\p$-linear automorphism $\vf: V\to V$ such that $\vf(L_1)=L_2$.
\medskip
Numbers $n$, $r$ are called the dimension and the rank of $L$
respectively.  Let $\{e_*\}=e_1, ..., e_r$ be a $\n F_q[\theta]$-basis of
$L$ such that $e_1, ..., e_n$ form a $\p$-basis of $V$.
\medskip
{\bf Definition 1.5.} The Siegel matrix of $L$ with respect to $\{e_*\}$ is the matrix $S=(s_{ij})\in M_{(r-n)\times n}(\p)$,  whose lines are coordinates of $e_{n+1}, ..., e_r$ in the basis $e_1, ..., e_n$:
$$\forall i =1,..., r-n \ \ \ \  e_{n+i}=\sum_{j=1}^n s_{ij}e_j\eqno{(1.5.1)}$$

{\bf 1.6.} We shall use the following convention for the action of $GL_{2n}(\n F_{q}[\theta])$ on the set of Siegel matrices. Let $L_1$, $L_2$ be lattices of rank $2n$ in $n$-dimensional vector spaces $V_1$, $V_2$ respectively, $\vf: V_1\to V_2$ a linear map such that $\vf(L_1)=L_2$. Let $g_1, \dots, g_{2n}$, $h_1, \dots, h_{2n}$ be $\n F_q[\theta]$-bases of $L_1$, $L_2$ respectively. We denote the matrix columns $(g_1,\dots, g_{2n})^t$, $(h_1,\dots, h_{2n})^t$ by $g_*$, $h_*$ respectively.  There exists a (uniquely defined) matrix $Z\in GL_{2n}(\n F_q[\theta])$ such that $$\vf(g_*)=Z^t h_*\eqno{(1.6.1)}$$

We denote by $\goth g_i$, $\goth h_i$, $\goth Z_{ij}$, $i, j=1,2$, the $i$-th (for $g$, $h$) and the $(i,j)$-th $(n\times 1)$, resp. $(n\times n)$-block of $g_*$, $h_*$, $Z$ respectively. This means that (1.6.1) becomes ($i=1,2$)
$$\vf(\goth g_i)=\goth Z_{1i}^t\goth h_1+\goth Z_{2i}^t\goth h_2\eqno{(1.6.2)}$$

We use (1.6.1), (1.6.2) in order to define the action of $GL_{2n}(\n F_q[\theta])$ on the set of Siegel matrices. Namely, let $S_1$, $S_2$ be the Siegel matrices of  $g_*$, $h_*$ respectively, i.e. $\goth g_2=S_1\goth g_1$, $\goth h_2=S_2 \goth h_1$. We let $S_2=Z(S_1)$. (1.6.2) implies the explicit formula
$$S_1=(\goth Z_{12}^t+\goth Z_{22}^tS_2)(\goth Z_{11}^t+\goth Z_{21}^tS_2)^{-1}\eqno{(1.6.3)}$$

{\bf Remark.} $S_2$ can be obtained as a function of $S_1$ by means of $\goth Z^{-1}$. We apologise for using of $\goth Z^t$ instead of $\goth Z$ itself (these notations appear because of duality between vectors and their coordinates ( = linear forms)).
\medskip
{\bf (1.6.4).} If $g_1, \dots, g_{2n}$ and $Z$ are given then $(h_1, \dots, h_{2n})$ is uniquely defined by (1.6.1). It can happen that $h_1, \dots, h_{n}$ is not a $\p$-basis of $V$. This is a condition $|\goth Z_{11}^t+\goth Z_{21}^tS_2|=0$, hence in this case the action of $Z$ on $S$ is not defined, i.e. we have only an "almost action" of $GL_{2n}(\n F_q[\theta])$ on the set of Siegel matrices. It is easy to see that this can happen even for $S=\om I_n$. We shall neglect this phenomenon, in all cases that we shall consider it does not exist, i.e. $Z(S)$ is defined.
\medskip
We have a result:
\medskip
(a) A Siegel matrix $S$ defines its lattice (denoted by $\goth L(S)$ ) uniquely (not all matrices in $ M_{(r-n)\times n}(\p)$ are Siegel matrices of lattices);
\medskip
(b) Two Siegel matrices $S_1$, $S_2$ define isomorphic lattices (i.e. $\goth L(S_1)=\goth L(S_2)$ ) if and only if there exists $Z\in GL_{2n}(\n F_q[\theta])$ such that $S_2=Z(S_1)$.
\medskip
{\bf 1.7.} Let us give more details on Remark 0.9.B. First, we have
\medskip
{\bf Proposition 1.7.1.} For any neighborhood $U$ of $\omega I_n$ in the set of Siegel matrices there exist $S_1$, $S_2\in U$, $\gamma\in GL_{2n}(\Bbb F_q[\theta])$ such that $S_2=\gamma(S_1)$ and $\gamma\not\in G_0$.
\medskip
{\bf Proof --- Example.} (A. Zobnin)\footnotemark \footnotetext{Found by a computer search, verified by hand calculation.}. We consider the case $q=2$, in this case $\omega\in \n F_4$ satisfies $\omega^2=\omega+1$ and
$G_0=\{\left(\matrix A&B\\ B&A+B \endmatrix \right)\}$. For simplicity, we consider the case $n=2$. Here and in (1.7.2), (1.7.3) we use another action of $GL_{2n}(\Bbb F_q[\theta])$ (not the one of (1.6.3)), namely $\ga(S):=(C+DS)(A+BS)^{-1}$. For any $m>0$ we let

$$\ga=\ga(m):=\left(\matrix \th^m+1&\th^m &&0&\th^m\\ \th^m+1&0&&\th^m&\th^m+1\\ \\ 0&\th^m+1&&\th^m+1&1\\ \th^m&\th^m+1&&0&\th^m+1\endmatrix \right)\eqno{(1.7.1.1)}$$ (spaces between its rows and columns indicate its $2\times2$-block structure),

$$S_1=S_2=S=S(m):=\left(\matrix \om+\th^{-m}&0 \\ 0&\om\endmatrix \right)\eqno{(1.7.1.2)}$$
\medskip
We have $\ga(S)=S$, $|\ga|=1$, $\ga\not\in G_0$, and ord $S(m)-\om I_2=m$. Finding of analogous examples for all $q$ and $n\ge2$ is an exercise for the reader. $\square$
\medskip
{\bf Remark.} For the present example the lattice defined by $S$ is reducible, moreover, it is a direct sum of two isogenous lattices of dimension 1 and rank 2. It is clear that a small deformation of $S$, with the same $\ga$, gives us a counterexample with an irreducible lattice.
\medskip
We see that the situation for the functional field case is not the same as for the number field case. The next proposition shows that not all is too bad:
\medskip
{\bf  Proposition 1.7.2.} There exists a neighborhood $U$ of $\omega I_n$ such that if $\ga\in GL_{2n}(\Bbb F_q[\theta])$ satisfies $\ga(\omega I_n)\in U$ then $\ga\in G_0$.
\medskip
{\bf  Proof.} We can choose $U=\{X|\ord(X-\omega I_n)>0\}$. Let $\ga= \left(\matrix A&B\\ C&D\endmatrix\right)$ ($n\times n$-blocks). Since $\ga(\omega I_n)\in M_n(\n F_{q^2}(\theta))$, we can denote $\ga(\omega I_n)=Y+(I_n+Z)\omega$, where $Y$, $Z\in M_n(\n F_{q}(\theta))$,  hence (here $q$ is odd, and $k:=\omega^2\in \n F_{q}$):
$$C+D\omega=(Y+(I_n+Z)\omega)(A+B\omega)$$ that implies
$$C=YA+k(I_n+Z)B\eqno{(1.7.2.1)}$$
$$D=YB+(I_n+Z)A\eqno{(1.7.2.2)}$$
This implies
$$\left(\matrix Y+(I_n+Z)\omega &-I_n\\ Y-(I_n+Z)\omega &-I_n\endmatrix\right) \left(\matrix A&B\\ C&D\endmatrix\right) \left(\matrix I_n&I_n\\ \omega I_n& -\omega I_n\endmatrix\right) = $$ $$\left(\matrix 0& 2\omega(I_n+Z)(A-B\omega)\\ -2\omega(I_n+Z)(A+B\omega)&0\endmatrix\right)$$
Condition $\ord Y, \ \ord Z>0$ implies $\ord \ \left|\matrix Y+(I_n+Z)\omega &-I_n\\ Y-(I_n+Z)\omega &-I_n\endmatrix\right|=0$. We have $|\ga|\in \n F_q^*$, hence $\ord |\omega(I_n+Z)(A+B\omega)|=0$, i.e. $\ord |A+B\omega|=0$. We have $|A+B\omega|\in \n F_{q^2}[\theta]$, hence $|A+B\omega|\in \n F_{q^2}^*$ and $(A+B\omega)^{-1}\in M_n(\n F_{q^2}[\theta])$. This means that $\ga(\om I_n)\in M_n(\n F_{q^2}[\theta])$ and hence, because $\ord (\ga(\om I_n) - \om I_n)>0$, we get $\ga\in G_0$. $\square$
\medskip
{\bf Remark 1.7.3.} There exists another proof of Proposition 1.7.2 for the symplectic group $GSp_{2n}(\Bbb F_q[\theta])$ (and symmetric Siegel matrices), this is the case of negatively self-dual Anderson t-motives, see [L], Section 7. Let us give it for completeness. We use the same notations.
\medskip
{\bf Proof for the symplectic case.} Multiplying (1.7.2.1) by $-B^t$ from the right and (1.7.2.2) by $A^t$ from the right and adding we get $$\la=(I_n+Z)(AA^t-kBB^t)$$ where $\la=DA^t-CB^t\in \n F_q^*I_n$.

Let $\al$ be an entry of $AA^t-kBB^t$ with the minimal ord, let it be the $(i,j)$-th entry, and let ord $(\al)=\de$. Since ord $Z>0$, we have ord $Z(AA^t-kBB^t)>\de$ and the $(i,j)$-th entry of $(I_n+Z)(AA^t-kBB^t)$ has ord $=\de$. Hence, $\de=0$, i.e. $AA^t-kBB^t\in M_n(\n F_q)$. We have $\la-(AA^t-kBB^t)=Z(AA^t-kBB^t)\in M_n(\n F_q)$, hence the conditions ord $Z>0$, ord $(AA^t-kBB^t)=0$ imply $Z(AA^t-kBB^t)=0$. This implies $\la=AA^t-kBB^t$ and hence $Z=0$.

Multiplying (1.7.2.1) by $A^t$ from the right and (1.7.2.2) by $-kB^t$ from the right and adding we get $CA^t-kDB^t=\la Y$. Since ord $Y>0$ and $CA^t-kDB^t\in M_n(\n F_q[\theta])$ this implies $Y=0$. $\square$

\medskip

{\bf 2. From a matrix $A$ to a Siegel matrix. }

\nopagebreak
\medskip
Recall that we consider $M$ given by the equation (0.8), $q$ is odd, and $\omega \in \n F_{q^2} - \n F_q$ satisfies $\omega^2\in \n F_q$. A Siegel matrix of $M(0)$ is $\omega I_n$. We consider 4 sets $\goth S_1, ... , \goth S_4$:
\medskip
$\goth S_1$. The set of $n \times n$ matrices $A$.
\medskip
$\goth S_2$. The set of t-motives $M$ given by the equation (0.8).
\medskip
$\goth S_3$. The set of $n \times n$ Siegel matrices $S$.
\medskip
$\goth S_4$. The set of lattices of rank $r=2n$ in $\p^n$.
\medskip
Let $\goth W_1\subset \goth S_1:=\{ A | \hbox{ord }A> \frac{q}{q^2-1}\}$, $\goth W_3\subset \goth S_3:=\{S | \hbox{ord }(S-\omega I_n)> \frac{q}{q^2-1}\}$ be open neighborhoods of 0, resp. $\omega I_n$ in $\goth S_1$, resp. $\goth S_3$. Let $M: \goth S_1 \to \goth S_2$, $\goth L: \goth S_3 \to \goth S_4$, $L:$ (a subset of $\goth S_2$ corresponding to uniformizable t-motives) $\to \goth S_4$ be as above. We shall show that all t-motives in $M(\goth W_1)$ are uniformizable, hence $L: M(\goth W_1)\to \goth S_4$ is defined. We have a diagram:
$$\matrix \goth S_1 & \hookleftarrow & \goth W_1 &\goth W_3 & \hookrightarrow & \goth S_3 \\ \\   M\downarrow & &  M\downarrow &&& \goth L\downarrow \\ \\  \goth S_2 & \hookleftarrow & M(\goth W_1) &\overset{L}\to{\to}& & \goth S_4
\endmatrix $$
\medskip
{\bf Proposition 2.}  There exists an isomorphism $\Cal S:\goth W_1\to\goth W_3 $ defined by power series of the entries of $A$ making the above diagram commutative.
\medskip
{\bf Proof.} As usual, for a matrix $A=(a_{ij})$ we denote $A^{(k)}=(a_{ij}^{q^k})$, and we denote $\theta_{ij}=\theta^{q^i}-\theta^{q^j}$.

For a given $A\in \goth S_1$ we denote the exponential map of $M(A)$ by $\Exp_A$. We have $\Exp_A(X)=\sum_{i=0}^\infty C_iX^{(i)}$ where $C_i=C_i(A)$, $C_0(A)=1$, and they satisfy the following recurrence relation (here $i\ge 1$; $C_{-1}(A)=0$):

$$C_\nu=\frac{AC_{\nu-1}^{(1)}+C_{\nu-2}^{(2)}}{\theta_{\nu0}}\eqno{(2.1)}$$
Let $J=(j_1,\dots,j_l)$ be a sequence of numbers $j_i\ge0$ or $J=\emptyset$. We denote $A^{(J)}:=A^{(j_1)}\cdot A^{(j_2)}\cdot\dots\cdot A^{(j_l)}$, $A^{(\emptyset)}=I_n$. We denote $m(J):=\max(l,j_1,\dots,j_l)$, $m(\emptyset):=0$. $J$ is called $\nu$-special if $J=\emptyset$ or if $j_1,\dots,j_l$ satisfy $0\le j_1<j_2<...<j_l<\nu$, and $J$ is called special if it is $\nu$-special for some $\nu$.  If $J$ is $\nu$-special then $m(J)\le \nu$, $m(J)\le j_l+1$.

We have (this follows from (2.1) immediately by induction): $C_\nu$ is a finite sum of terms of the form
$$\frac{1}{\theta_{\nu,i_1}\cdot\theta_{\nu,i_2}\cdot\dots\cdot\theta_{\nu,i_k}\cdot\theta_{\nu0}} A^{(J)}\eqno{(2.2)}$$
where $\nu>i_1>i_2>... >i_k>0$ is a sequence of numbers uniquely defined by $J$, and $J$ is $\nu$-special. We denote $\xi_{\nu,J}:=\frac{1}{\theta_{\nu,i_1}\cdot\theta_{\nu,i_2}\cdot\dots\cdot\theta_{\nu,i_k}\cdot\theta_{\nu0}}$, hence (2.2) becomes $$C_\nu(A)=\sum_J \xi_{\nu,J}A^{(J)}\eqno{(2.3)}$$ the sum runs over $\nu$-special $J$. For any fixed $J$ there is no more than one term of type (2.2) in $C_\nu$, and $$k+1\ge \lceil \frac{\nu}2\rceil\eqno{(2.4)}$$ where $\lceil x \rceil:=\min \{\al\in \n Z|\al\ge x\}$ is the ceiling function.
Further, the only terms corresponding to $J=\emptyset$ are
$$\xi_{\nu,\emptyset}=\frac{1}{\theta_{\nu,\nu-2}\cdot\theta_{\nu,\nu-4}\cdot\dots\cdot\theta_{\nu2}\cdot\theta_{\nu0}} \eqno{(2.5)}$$ for even $\nu$, the only terms corresponding to the case $l=1$, $j_1=0$ (i.e. $J=(0)$ ) are terms
$$\xi_{\nu,(0)}A\hbox{ where }\xi_{\nu,(0)}=\frac{1}{\theta_{\nu,\nu-2}\cdot\theta_{\nu,\nu-4}\cdot \dots\cdot\theta_{\nu3}\cdot\theta_{\nu1}\cdot\theta_{\nu0}}\eqno{(2.6)}$$ for odd $\nu$ ($\xi_{\nu,(0)}=0$ for even $\nu$).

Recall that the exponent for the Carlitz module $\goth C_2$ has the form
$$\Exp_0(z)=z+\frac{1}{\theta_{20}}z^{q^2}+\frac{1}{\theta_{42}\theta_{40}}z^{q^4}+...=\sum_{i=0}^\infty\  \xi_{2i,\emptyset} \ z^{q^{2i}}\eqno{(2.7)}$$
We denote by
$y_0\in \p$ a nearest-to-zero root to $\Exp_0(z)=0$ (this is $\xi$ of $\goth C_2$ in notations of [G]). It is defined up to
multiplication by elements of $\n F_{q^2}^*$, and it generates over $\n F_{q^2}[\theta]$ the
lattice of the Carlitz module $\goth C_2$. We fix one such $y_0$. We have ord$(y_0)=-\frac{q^2}{q^2-1}$, it corresponds to the first segment of the Newton polygon of (2.7).

We identify $\Lie(\goth C_2)$ (see [G], one line above Definition 5.4.5 for the Lie space of a t-motive) with $\p$ and hence $\Lie(M(0))$ with $\p^n$. We denote the standard $\p$-basis of $\Lie(M(0))$ by $l_{1}, ... , l_{n}$, namely $l_{i}=(0,...,0,1,0,...,0)$ (1 at the $i$-th place),
we denote the elements $y_0l_i$ by $\goth e_i$, and let $\goth e$ be the column $(\goth e_1,\dots,\goth e_n)^t$. Also for $i=1,\dots, n$ we denote $g_i(0):=\goth e_i$, $g_{n+i}(0):=\omega \goth e_i$, and we denote by $Y_i$, resp. $Y'_i$ the vector columns of coordinates of $g_i(0)$, resp. $g_{n+i}(0)$, i.e. $Y_{i}=(0,...,0,y_0,0,...,0)^t$, resp. $Y'_{i}=(0,...,0,\omega y_0,0,...,0)^t$ ($y_0$, resp. $\omega y_0$ at the $i$-th place). Hence, elements $g_i(0)$ for $i=1,\dots, 2n$ form a basis of $L(M(0))$ over $\n F_q[\theta]$. Further, like in (1.6), we denote the column $(g_1(0), \dots, g_n(0))^t$, resp. the column $(g_{n+1}(0), \dots, g_{2n}(0))^t$ by $\goth g_1(0)$, resp. $\goth g_2(0)$.

Now let us consider the deformation of this basis for $A\in \goth W_1$. We let  $$X_i=Y_i+\Delta_i, \ \  \ X'_i=Y'_i+\Delta'_i\eqno{(2.8)}$$ where $X_i$, $X'_i$, $\Delta_i$, $\Delta'_i$ depend on $A$ and where $X_i$, $X'_i$ are solutions to the equation $$\sum_{j=0}^\infty C_j(A) X^{(j)}=0\eqno{(2.9)}$$ near $Y_i$, $Y'_i$ respectively (i.e. $\Delta_i$, $\Delta'_i$ are small complements). We shall show that they exist and are unique. We denote vectors, whose coordinates are columns $X_i$, $X'_i$, by $g_i(A)$, $g_{n+i}(A)$ respectively, and we denote the column $(g_1(A), \dots, g_n(A))^t$, resp. the column $(g_{n+1}(A), \dots, g_{2n}(A))^t$ by $\goth g_1(A)$, resp. $\goth g_2(A)$.
This will mean that we get a lattice generated by $g_i(A)$, $g_{n+i}(A)$. The Siegel matrix of the basis $g_i(A)$, $g_{n+i}(A)$ is, by definition, $\Cal S(A)$.

To prove that $\Cal S$ is an isomorphism we need to find an explicit expression for $\Cal S(A)$. Let us denote $D_j=D_j(A)=C_j(A)-C_{j}(0)$. Substituting (2.8) to (2.9) we get $$\sum_{j=1}^\infty D_j Y_i^{(j)} +\Delta_i +\sum_{j=1}^\infty C_j\Delta_i^{(j)}=0\eqno{(2.10)}$$
$$\sum_{j=1}^\infty D_j {Y'_i}^{(j)} +\Delta'_i +\sum_{j=1}^\infty C_j{\Delta'_i}^{(j)}=0\eqno{(2.11)}$$
where $C_j=C_j(A)$, $D_j=D_j(A)$ are parameters, $\Delta_i$, $\Delta'_i$ are column matrix unknowns.
We shall see later that for small $A$ the sums $\sum_{j=1}^\infty D_j Y_i^{(j)}$, $\sum_{j=1}^\infty D_j {Y'_i}^{(j)}$ converge, and both (2.10) and (2.11) have a unique solution near 0.

Now let us consider the $n\times n$-matrix form of (2.10), (2.11). Let $Y$, resp. $Y'$, $\Delta$, $\Delta'$ be $n\times n$-matrices whose $i$-th column is $Y_i$, resp. $Y_i'$, $\Delta_i$, $\Delta'_i$. We have $Y=y_0I_n$, $Y'=\omega y_0I_n$. The $n\times n$-matrix form of (2.10), (2.11) is the following:
$$\sum_{j=1}^\infty D_j Y^{(j)} +\Delta +\sum_{j=1}^\infty C_j\Delta^{(j)}=0\eqno{(2.12)}$$
$$\sum_{j=1}^\infty D_j {Y'}^{(j)} +\Delta' +\sum_{j=1}^\infty C_j{\Delta'}^{(j)}=0\eqno{(2.13)}$$

We need a high-dimensional version of the classical Hensel lemma:
\medskip
{\bf Lemma 2.14.} Let $U_0+X+U_1X^{(1)}+U_2X^{(2)}+\dots=0$ be a matrix equation (i.e. $X\in M_n(\p)$ is an unknown matrix, $U_i\in M_n(\p)$ are matrix parameters) such that ord $U_0>0$, ord $U_i\ge0$. Then this equation has a unique solution $X$ satisfying ord $X>0$.
\medskip
{\bf Proof} of existence is completely analogous to the proof of the classical Hensel lemma. Unicity is also obvious: if we have 2 solutions $X$, $X'$ then $-X+X'=\sum_{i=1}^\infty U_i(X-X')^{(i)}$. Comparing the maximal ord of entries of the left and right hand sides of this equality we get immediately $X-X'=0$. $\square$
\medskip
We shall need the explicit form of the solution to the equation of Lemma 2.14. It is $$X=-U_0+\sum \pm U_{\alpha_1}^{(\beta_1)}\cdot U_{\alpha_2}^{(\beta_2)}\cdot\dots \cdot U_{\alpha_{\la-1}}^{(\beta_{\la-1})}\cdot U_0^{(\beta_\la)}\eqno{(2.15)}$$ where $\la\ge2$, $\alpha_i>0$ for $i=1,\dots,\la-1$, $\be_\la\ge1$, $\be_i\ge0$ for $i=1,\dots,\la-1$, and for any given $\nu$ there is only finitely many terms such that $\beta_\la\le \nu$. This follows immediately by induction applied to the process of solution of 2.14.
\medskip
Let us evaluate the order of coefficients. From here until 2.19 we consider only special $J$. For (2.12) we have: $U_0=\sum_{j=1}^\infty D_j Y^{(j)}$, and (2.3) implies $U_0=\sum_{J\ne\emptyset}(\sum_{\nu\ge0} \xi_{\nu,J}Y^{(\nu)})A^{(J)}$. Really, in the inner sum we have $\nu\ge m(J)$. (2.4) implies that ord $\xi_{\nu,J}\ge \lceil\frac{\nu}2\rceil q^\nu$, hence ord $\xi_{\nu,J}Y^{(\nu)}\ge q^\nu(\lceil\frac{\nu}2\rceil  - \frac{q^2}{q^2-1})$ which is $\ge q^\nu$ for $\nu\ge3$. This means that for any $J$ the sum $\bar \xi_{0,J}:=\sum_\nu \xi_{\nu,J}Y^{(\nu)}$ converges, $$U_0=\sum_{J\ne\emptyset}\bar \xi_{0,J}A^{(J)}\eqno{(2.16)}$$ and if $m(J)\ge3$ then ord $\bar \xi_{0,J}\ge q^{m(J)}$. If $m(J)\le2$ then all terms of the sum $\sum_\nu \xi_{\nu,J}Y^{(\nu)}$ having $\nu\ge3$ have ord $\ge q^3\frac{q^2-2}{q^2-1}$. The only term having $\nu=2$ is the term $\frac{1}{\theta_{21}\theta_{20}}Y^{(2)}$ for $J=(0,1)$. Its ord is $\ge q^2\frac{q^2-2}{q^2-1}$. The only term having $\nu=1$ corresponds to $J=(0)$. Hence, for all $J\ne (0)$ we have: ord $\bar \xi_{0,J}\ge q^{m(J)}\frac{q^2-2}{q^2-1}$.

For $J=(0)$ we have (see (2.6)) $\bar \xi_{0,(0)}=d$ where
$$d=\frac{y_{0}^{q}}{\theta_{10}}+\frac{y_{0}^{q^3}}{\theta_{31}\theta_{30}}
+ \frac{y_{0}^{q^5}}{\theta_{53}\theta_{51}\theta_{50}}+
\frac{y_{0}^{q^7}}{\theta_{75}\theta_{73}\theta_{71}\theta_{70}}+...=\sum_{i=1}^\infty \xi_{i,(0)}y_{0}^{q^i}$$
Later we shall use that $\bar \xi_{0,(1)}=0$ --- really, (2.1) obviously implies that $\forall \nu$ we have $\xi_{\nu,(1)}=0$.

We have ord $d=\ord \frac{1}{\theta_{10}}Y^{(1)}=-\frac{q}{q^2-1}$. We fix a small $\ve>0$ and we change the scale, letting $A=\mu B$ where $\mu$ satisfies ord $\mu=\frac{q}{q^2-1}+\ve$. We define $\hat \xi_{0,J}$ by the substitution $A=\mu B$ to (2.16) to get
$$U_0=\sum_{J\ne\emptyset,(1)} \hat \xi_{0,J} B^{(J)}\eqno{(2.16.1)}$$ we have $\hat \xi_{0,J}=\bar \xi_{0,J}\ \mu^{\ga}$ for some integer $\ga>0$. Obviously $$\ord \hat \xi_{0,J}\ge \ve m(J)\eqno{(2.17)}$$ (we do not need more strong inequality). We see that if ord $B\ge0$ then (2.16.1) converges.

We define numbers $\hat \xi_{\nu,J}$ by the equality $C_\nu=\sum_J \hat \xi_{\nu,J} B^{(J)}$ obtained by substitution of $A=\mu B$ to (2.3) (here $J$ can be $\emptyset$ or $(0)$). As earlier $$\ord \hat \xi_{\nu,J}\ge \ve m(J)\eqno{(2.18)}$$

{\bf (2.19)} To apply (2.15) we need more notations (from here $J$ is not necessarily special). We define $(J+\beta):=(j_1+\beta,\dots,j_m+\beta)$, and if $(J_1)=(j_{11},\dots,j_{1,m_1})$, $(J_2)=(j_{21},\dots,j_{2,m_2})$, then $(J_1\cup J_2):=(j_{11},\dots,j_{1,m_1},j_{21},\dots,j_{2,m_2})$. Applying (2.15) to (2.12) we get a formal sum
$$\Delta=-d\mu B-\sum_{J\ne \emptyset,(0),(1)} \hat \xi_{0,J}  B^{(J)}+\sum_{\la,J_*,\alpha_*,\beta_*}\pm\hat \xi_{\alpha_*,J_*}^{(\beta_*)} B^{(J_*+\beta_*)}\eqno{(2.20)}$$ where $J_*=(J_1,\dots,J_\la)$ and $J_i=(j_{i1},\dots,j_{i,m_i})$, $\alpha_*=(\alpha_1, \dots,\alpha_\la)$ and $\alpha_\la=0$, $\beta_*=(\beta_1, \dots, \beta_\la)$, $\hat \xi_{\alpha_*,J_*}^{(\beta_*)}:=\prod_{i=1}^\la \hat \xi_{\alpha_i,J_i}^{(\beta_i)}$, $(J_*+\beta_*):=(J_1+\beta_1)\cup\dots\cup(J_\la+\beta_\la)$.
\medskip
There is no term $\xi_*B^{(1)}$ in (2.20). Really, there is no such term in $U_0$. If a term of the sum of (2.15) has $\be_\la\ge2$ then it cannot contain a term $k_*B^{(1)}$. The process of solution of the equation of (2.14) shows that a unique term of the sum of (2.15) having $\be_\la=1$ is the term $U_1U_0^{(1)}$. Since $U_1=\frac1{\theta_{10}}A$ we get that it does not contain a term $\xi_*B^{(1)}$ as well.

Let us evaluate ord of coefficients. We need two elementary lemmas. For a term $k B^{(J)}$ condition (*) means ord $k\ge \ve m(J)$.
\medskip
{\bf Lemma 2.21.} Let a term $\xi B^{(J)}$ satisfies (*). Then the term $(\xi B^{(J)})^{(\beta)}$ satisfies (*).
\medskip
{\bf Proof.} $(\xi B^{(J)})^{(\beta)}=\xi^{q^\beta}B^{(J+\beta)}$. We have $m(J)+\beta\ge m(J+\beta)$, $q^\beta \cdot m(J)\ge m(J)+\beta$ (if $m(J)\ge1$), hence ord $\xi^{q^\beta}=q^\beta\ord \xi\ge q^\beta\ve m(J)\ge \ve(m(J)+\beta)\ge\ve m(J+\beta)$. If $m(J)=0$, i.e. $J=\emptyset$, the lemma obviously holds. $\square$
\medskip
{\bf Lemma 2.22.} Let terms $\xi_1B^{(J_1)}$, $\xi_2 B^{(J_2)}$ satisfy (*). Then their product also satisfies (*).
\medskip
{\bf Proof.} We have $m(J_1)+m(J_2)\ge m(J_1\cup J_2)$ (this follows from $\max (m_1, \ga_1)+\max (m_2, \ga_2)\ge \max(m_1+m_2, \max(\ga_1,\ga_2))$ if all entries are $\ge0$; this inequality follows from
\medskip
$\max (m_1, \ga_1)+\max (m_2, \ga_2)\ge m_1+m_2$,

$\max (m_1, \ga_1)+\max (m_2, \ga_2)\ge \ga_1$,

$\max (m_1, \ga_1)+\max (m_2, \ga_2)\ge \ga_2$).

Adding  ord $\xi_1\ge \ve m(J_1)$ and ord $\xi_2\ge \ve m(J_2)$ we get the result. $\square$
\medskip
Now let us prove that for any $\nu$ there is only finitely many terms of (2.20) such that for its $J$ we have $m(J)\le \nu$. If $m(J)\le \nu$ then $\be_\la$ of (2.15) is $\le \nu$. We have  only finitely many terms of (2.15) satisfying $\be_\la\le \nu$. For any of these terms the factors $U_{\alpha_i}^{(\beta_i)}$ for $i<\la$ contain only finitely many terms, and the $\la$-th factor $U_{0}^{(\beta_\la)}$ contains only finitely many terms having $m(J)\le \nu$, because of (2.16.1).
\medskip
Further, (2.17) and (2.18) imply that for all $J\ne (0)$, for all $\nu=0,1,\dots$ we have $\ord \hat \xi_{\nu,J}\ge\ve m(J)$. Since $\ve = \ord d\mu$, the same holds for $J=(0)$ as well. The above lemmas show that if we define $\goth k_{J}$ writing (2.20) in the form $$\Delta=-d\mu B+\sum_{J\ne\emptyset,(0), (1)} \goth k_J  B^{(J)} \eqno{(2.23)}$$ then ord $\goth k_J\ge \ve m(J)$. The same formulas hold for $\Delta'$: $$\Delta'=\omega d\mu B+\sum_{J\ne\emptyset,(0), (1)} \goth k'_J  B^{(J)}\eqno{(2.24)}$$  (it is clear that the coefficient $\goth k'_{(0)}$ is really $\omega d\mu$, because $\sum_{i=1}^\infty \xi_{i,(0)}(\omega y_{0})^{q^i}=-\omega d$,  because $\omega^{q^i}=-\omega$ for odd $i$) and ord $\goth k'_J\ge \ve m(J)$.

(1.6) implies $S=\Cal S(A)=(\omega y_0I_n+{\Delta'}^t)(y_0I_n+\Delta^t)^{-1}=(\omega I_n+y_0^{-1}{\Delta'}^t)(I_n+y_0^{-1}\Delta^t)^{-1}$. Since if ord $B>0$ then ord $\Delta\ge \ve$, we have that
$$\frac{y_0 (S^t-\omega I_n)}{2\omega d\mu}=B - \frac1{2d\mu }\sum_{J\ne\emptyset,(0), (1)} \goth k_J  B^{(J)} + \frac1{2\omega d\mu } \sum_{J\ne\emptyset,(0), (1)} \goth k'_J  B^{(J)}$$ $$+  \frac1{2 d\mu } \sum_{i=2}^\infty \pm y_0^{-(i-1)}\Delta^i+\frac1{2\omega d\mu } \sum_{i=2}^\infty \pm y_0^{-(i-1)}\Delta^{i-1}\Delta' \eqno{(2.25)}$$ which can be written as $$\frac{y_0 (S^t-\omega I_n)}{2\omega d\mu}=B+\sum_{J\ne\emptyset,(0), (1)} K_J  B^{(J)}\eqno{(2.26)}$$ where $K_J$ are some coefficients. Really, it is clear that for any $J$ the sums in (2.25) contain only finitely many terms having $B^{(J)}$, because if $m(J)=m$ then $i$ in (2.25) is $\le m$. Further, Lemma 2.22 implies that ord $K_J\ge \ve (m(J)-1)$. Since the sum runs over $J$ having $m(J)\ge2$ we get that for $\ve_1:=\ve/2$ we have $$\ord K_J\ge \ve_1 \ m(J)\eqno{(2.27)}$$

(2.26) shows that $S$ is a function in $B$: $S=f(B)$. The below version of the Hensel lemma will show us that $f$ is 1 -- 1 in a neighborhood of 0. We need more general power series. Let $J=(j_1,\dots,j_\mu)$ be as above, it is called the type of a term of a series, and let $i=1,2,...$ be an integer parameter. Let $u_{iJ0}, \dots, u_{iJ\mu}\in M_n(\p)$ be coefficients associated to $J$ and $i$. The $i$-th term of type $J$ with coefficients $u_{iJ*}$ is, by definition, $$u_{iJ0}X^{(j_1)} u_{iJ1}X^{(j_2)}\cdot \dots \cdot u_{iJ,\mu-1}X^{(j_\mu)} u_{iJ\mu}$$ where $X\in M_n(\p)$ is a variable. We denote it by $C(J,i)(X)$ (coefficients $u_{iJ*}$ are by default). The unique term ($i=1)$ corresponding to $J=\emptyset$ we denote simply by $u$. Finally, we denote $\sum_{\alpha=0}^\mu \ord u_{iJ\alpha}$ by $\ord (C(J,i))$.

Let us consider the power series $\sum_{J,i} C(J,i)(X)$ and the equation $$X=u+\sum_{J\ne\emptyset}\sum_i C(J,i)(X)\eqno{(2.28)}$$ (the sum can contain terms having $J=(0)$ ).
\medskip
{\bf Lemma 2.29.} Let there exist $\gamma_\emptyset, \ \gamma>0$ such that the following conditions hold:
\medskip
1. ord $u>\gamma_\emptyset$.
\medskip
2. $\forall J$, $\forall i$ we have $\ord (C(J,i))\ge\ga\ m(J)$.
\medskip
3. For all fixed $J$ we have: $\ord (C(J,i))\to +\infty$ as $i\to +\infty$.
\medskip
Then $\sum_{J,i} C(J,i)(X)$ converges if $\ord X\ge0$, and the equation (2.28) has a unique solution satisfying $\ord X\ge\gamma_\emptyset$.
\medskip
{\bf Proof. } We shall show that applying a step of successive approximation we come to an equation of the same type, satisfying to the same conditions, and that this method converges. First, let us show that $\sum_{J,i} C(J,i)(X)$ converges if $\ord X\ge 0$. We have $$\ord (C(J,i)(X))\ge \ord (C(J,i))+m(J) \ \ord X$$ hence (3) implies that $\forall J$ \ $\sum_i C(J,i)(X)$ converges, and (2) implies that $\forall J$
$$\ord \sum_i C(J,i)(X)\ge m(J)(\ga+\ord X)\eqno{(2.29.1)}$$ Since for a given $m$ there exists only finitely many $J$ having $m(J)\le m$ we have that $\sum_J \sum_i C(J,i)(X)$ converges.

We shall denote by prime the objects obtained after iteration. For example, the new unknown will be denoted by $X'$ (it satisfies $X'=X-u$), the new coefficients will be denoted by $u_{iJl}'$ etc. In order to prove existence of the solution (i.e. that the iteration process converges) it is sufficient to prove:
\medskip
(a). The iterated equation has the same form (2.28) (clearly with other coefficients $u'_{iJl}$);
\medskip
(b). Properties (1) - (3) hold for the iterated equation for the following values of $\gamma_\emptyset'$, $\ga'$:
\medskip
(c) $\gamma_\emptyset'=\ga_\emptyset+\gamma$, $\gamma'=\gamma$.
\medskip
Really, a term $C(J,i)(X)$ after the substitution $X=u+X'$ becomes
\medskip
$$u_{iJ0}(u^{(j_1)}+{X'}^{(j_1)}) u_{iJ1}(u^{(j_2)}+{X'}^{(j_2)})\cdot \dots \cdot u_{iJ,\mu-1}(u^{(j_\mu)}+{X'}^{(j_\mu)} )u_{iJ\mu}=$$ $$\sum_{J_0\subset J}C'(J_0,i')(X')\eqno{(2.29.2)}$$ where the sum runs over all $2^\mu$ subsequences of the sequence $J$, and all these terms are of the same type, this proves (a).
\medskip
We have $u'=\sum_{J}\sum_iC(J,i)(u)$. According (2.29.1), we have $\ord u' \ge \ga +\ga_\emptyset= \gamma_\emptyset'$. For any $J_0\subset J$ and $C'(J_0,i')$ from (2.29.2) we have $\ord C'(J_0,i')\ge \ord C(J,i)\ge \ga m(J)\ge \ga\ m(J_0)$. This proves (2) for the iterated equation. Let us fix $N$ and $J_0$, and prove that there is only finitely many $i'$ such that $\ord C'(J_0,i')\le N$. Again because $\ord C'(J_0,i')\ge \ord C(J,i)$ we see that in this case $\ord C(J,i)\le N$. There is only finitely many such $J$, $i$, and each of them gives only finitely many  $ C'(J_0,i')$. This proves existence of the solution.

Unicity:  Let $X_1$, $X_2=X_1+D$ be solutions to (2.28), where ord $D\ge\ga_\emptyset$. This means that $D=\sum_{J\ne\emptyset}\sum_i [C(J,i) (X_1+D)- C(J,i)(X_1)]$. (2.29.2) shows that $C(J,i) (X_1+D)- C(J,i)(X_1)$ is a sum of $2^\mu-1$ terms. The ord of each of these terms is greater than ord $D$ --- a contradiction. $\square$
\medskip
We apply this lemma to (2.26). We have $u=\frac{y_0 (S^t-\omega I_n)}{2\omega d\mu}$ has ord $>0$, 2.29 (2) is 2.27, and 2.29 (3) is true, because (2.26) contains only one term $K_J  B^{(J)}$ for any $J$. All conditions of the lemma are satisfied, hence we have a 1 -- 1 map $B \mapsto S$ in a neighborhood of 0. This proves Proposition 2. $\square$
\medskip
{\bf 3. Some notations. }
\medskip
Let $G_1=GL_n(\n F_{q^2}[\theta])$. Recall that $T_0=\theta+\tau^2\in M_n(\p\{\tau\})$ and $G_2$ are defined in (0.9.3),  (0.9.5). It is easy to see that $G_2=GL_n(\n F_{q^2}[T_0])\subset M_n(\p\{\tau\})$. There exists an isomorphism $\alpha: G_1\to G_2$ defined by $\alpha(\theta)=T_0$. We need an explicit formula for $g\in GL_n(\n F_{q^2}[T_0])$ as an element of $GL_n(\p\{\tau\})$. Namely, we define $k_{ij}\in \n F_q[\theta]$ (generalized binomial coefficients, $i,j\ge0$, $i\le j$) as follows:
$$T_0^j=\sum_{i=0}^j k_{ij} \tau^{2i}$$
Numbers $k_{ij}$ are defined either by recurrent relation $k_{ij}=k_{i-1,j-1}+k_{i,j-1}\theta^{q^{2i}}$, $k_{ii}=1$ for $i\ge 0$, $k_{-1,j}=0$ for $j\ge0$, or by the explicit formula $$k_{ij}=\sum_{P\in \Cal S_{ij} }\theta^{P(q^2)}=\sum_{P\in \Cal S_{ij} }\theta^{\overset{i}\to{ \underset{\gamma=0}\to{\sum}}c_\gamma q^{2\gamma}}\eqno{(3.1)}$$
where $\Cal S_{ij}\subset \n Z^+[x]$ is the set of polynomials $P(x)=\overset{i}\to{ \underset{\gamma=0}\to{\sum}}c_\gamma x^\gamma$ of degree $\le i$ such that $c_\gamma\ge0$, $\overset{i}\to{ \underset{\gamma=0}\to{\sum}}c_\gamma=j-i$.
\medskip
For $\goth B_1=\sum_{i=0}^m B_{1i}\theta^i\in G_1$, where $B_{1i}\in M_n(\n F_{q^2})$, we have
$$\alpha(\goth B_1)=\sum_{i=0}^m B_{1i}T_0^i=\sum_{i=0}^m(\sum_{j=i}^m k_{ij}B_{1j})\tau^{2i}\eqno{(3.2)}$$
\medskip
Example: Table of $k_{ij}$:
\medskip
$\left(\matrix k_{00}& k_{01}& k_{02}& k_{03}& k_{04} \\k_{11}& k_{12}& k_{13}& k_{14}& k_{15} \\k_{22}& k_{23}& k_{24}& k_{25}& k_{26}  \\k_{33}& k_{34}& k_{35}& k_{36}& k_{37}\endmatrix\right)=$
\medskip
\medskip
\noindent
$\matrix 1&\theta&\theta^{2}&\theta^{3}&\theta^{4}\\ \\ 1&\theta^{q^2}+\theta &\theta^{2q^2}+\theta^{q^2+1}+\theta^{2}&\theta^{3q^2}+\theta^{2q^2+1}+\theta^{q^2+2}
+\theta^{3}&\dots\\  \\
1&\theta^{q^4}+\theta^{q^2}+\theta&\theta^{2q^4}+\theta^{q^4+q^2}+\theta^{q^4+1}+&\theta^{3q^4}+\dots&\dots
\\ &&+\theta^{2q^2}+\theta^{q^2+1}+\theta^{2}\\ \\ 1&\theta^{q^6}+\theta^{q^4}+\theta^{q^2}+\theta&\theta^{2q^6}+\dots&\theta^{3q^6}+\dots&\dots\endmatrix$
\medskip
{\bf (3.3)} We have ord $(k_{ij})=-(j-i)q^{2i}$.
\medskip
We need also numbers $l_{ij}$ defined by a formula similar to (3.1):
$$l_{ij}=\sum_{P\in \Cal S_{ij} }\theta^{c_0+\overset{i}\to{ \underset{\gamma=1}\to{\sum}}c_\gamma q^{2\gamma-1}}\eqno{(3.4)}$$
(notations are the same as in (3.1)).
\medskip
{\bf Example:} Table of $l_{ij}$:
\medskip
$\left(\matrix l_{00}& l_{01}& l_{02}& l_{03}& l_{04} \\l_{11}& l_{12}& l_{13}& l_{14}& l_{15} \\l_{22}& l_{23}& l_{24}& l_{25}& l_{26}  \\l_{33}& l_{34}& l_{35}& l_{36}& l_{37}\endmatrix\right)=$
\medskip
\medskip
\noindent
$\matrix 1&\theta&\theta^{2}&\theta^{3}&\theta^{4}\\ \\ 1&\theta^{q}+\theta &\theta^{2q}+\theta^{q+1}+\theta^{2}&\theta^{3q}+\theta^{2q+1}+\theta^{q+2}
+\theta^{3}&\dots\\  \\
1&\theta^{q^3}+\theta^{q}+\theta&\theta^{2q^3}+\theta^{q^3+q}+\theta^{q^3+1}+&\theta^{3q^3}+\dots&\dots
\\ &&+\theta^{2q}+\theta^{q+1}+\theta^{2}\\ \\ 1&\theta^{q^5}+\theta^{q^3}+\theta^{q}+\theta&\theta^{2q^5}+\dots&\theta^{3q^5}+\dots&\dots\endmatrix$
\medskip
{\bf (3.5)} We have for $i>0$: ord $(l_{ij})=-(j-i)q^{2i-1}$, ord $(l_{0j})=-j$.
\medskip
We shall need the following technical
\medskip
{\bf Lemma 3.6.} $\forall  i,j\ge0$, $i\le j$ we have $l_{ij}=k_{ij}^q-\theta_{2i+1,0}l_{i+1,j}$ (we let $l_{i+1,i}=0$).
\medskip
{\bf Proof.} We consider the following subsets\footnotemark \footnotetext{Notation $\Cal S_{ij}^{(*)}$ has nothing common with the elevation to $q^*$-th power.} of the sets $\Cal S_{ij}$:
\medskip
$\Cal S_{ij}^{(1)}:=\{P\in\Cal S_{ij} \ |\ c_0=0\}$.
\medskip
$\Cal S_{ij}^{(2)}:=\{P\in\Cal S_{ij}\ |\ c_0\ne0\}$. We have  \ $\Cal S_{ij}^{(1)}\cup\Cal S_{ij}^{(2)}=\Cal S_{ij}$, $\Cal S_{ij}^{(1)}\cap\Cal S_{ij}^{(2)}=\emptyset$.
\medskip
$\Cal S_{ij}^{(3)}:=\{P\in\Cal S_{ij} \ |\ c_i=0\}$.
\medskip
$\Cal S_{ij}^{(4)}:=\{P\in\Cal S_{ij}\ |\ c_i\ne0\}$. We have \ $\Cal S_{ij}^{(3)}\cup\Cal S_{ij}^{(4)}=\Cal S_{ij}$, $\Cal S_{ij}^{(3)}\cap\Cal S_{ij}^{(4)}=\emptyset$.

\medskip
For $P=\overset{\delta}\to{ \underset{\gamma=0}\to{\sum}}c_\gamma x^\gamma\in\Cal S_{\delta j}$ we denote $\goth M_k(P):=\theta^{\overset{\delta}\to{ \underset{\gamma=0}\to{\sum}}c_\gamma q^{2\gamma}}$, $\goth M_l(P):=\theta^{c_0+\overset{\delta}\to{ \underset{\gamma=1}\to{\sum}}c_\gamma q^{2\gamma-1}}$ --- the monomials corresponding to $P$ that enter in (3.1), (3.4) respectively (here $\delta=i$ or $i+1$). There are isomorphisms of sets:
\medskip
$\beta_1:\Cal S_{ij}^{(1)} \to \Cal S_{ij}^{(3)}$, $\beta_1(P):=P/x$;
\medskip
$\beta_2:\Cal S_{ij}^{(2)} \to \Cal S_{i+1,j}^{(3)}$, $\beta_2(P):=P-1$;
\medskip
$\beta_3:\Cal S_{ij}^{(4)} \to \Cal S_{i+1,j}^{(1)}$, $\beta_3(P):=Px-x^{i+1}$;
\medskip
$\beta_4:\Cal S_{i+1,j}^{(2)} \to \Cal S_{i+1,j}^{(4)}$, $\beta_4(P):=P+x^{i+1}-1$.
\medskip
The following formulas are checked immediately:
\medskip
$\forall \ P\in \Cal S_{ij}^{(1)}$ we have: $\goth M_l(P)=\goth M_k(\beta_1(P))^q$;
\medskip
$\forall \ P\in \Cal S_{ij}^{(2)}$ we have: $\goth M_l(P)=\theta \ \goth M_l(\beta_2(P))$;
\medskip
$\forall \ P\in \Cal S_{ij}^{(4)}$ we have: $\goth M_k(P)^q=\theta^{q^{2i+1}} \ \goth M_l(\beta_3(P))$;
\medskip
$\forall \ P\in \Cal S_{i+1,j}^{(2)}$ we have: $\theta^{q^{2i+1}} \ \goth M_l(P)=\theta \ \goth M_l(\beta_4(P))$.
\medskip
These formulas imply the lemma. $\square$
\medskip
{\bf 4. Main theorem. }
\medskip
We recall that $G_0\subset GL_{2n}(\n F_q[\theta])$ is the stabiliser of the Siegel matrix $\omega I_n$. 
(1.6.3) implies that $G_0\subset \{\ga=\left(\matrix C&\omega^2D\\ D&C \endmatrix\right)\}$ such that $C,D\in M_n(\n F_q[\theta])$. There exists an isomorphism $\beta: G_0 \to G_1=GL_n(\n F_{q^2}[\theta])$ defined by the formula $\beta(\gamma)=C+\omega D$.
\medskip
Let $\Delta_A$, $\Delta_X\in M_n(\p)$ be small. The action of the above $\gamma$ in a neighborhood of $\omega I_n$ is described by
\medskip
{\bf Lemma 4.1.} If $\omega+\Delta_X=\gamma(\omega+\Delta_A)$ then

$$\Delta_A=(C-\omega D)^t\Delta_X(C+\omega D)^{t-1}+\hbox{ higher terms in }\Delta_X  \eqno{\square}$$

Let $\ga$, $m$ be from (0.9.2). We let $\goth B_1:=\be(\ga)\in G_1$, $\goth B_2:=\al\circ\be(\ga)\in G_2$, $\goth B_1=\sum_{i=0}^m B_{1i}\theta^i\in G_1$, where $B_{1i}\in M_n(\n F_{q^2})$ like in (3.2). Now we can state the main calculational proposition. Let $W_0:=m^qq^{2m+1}$, $W:=(2m+2)W_0$.
\medskip
{\bf Proposition 4.2.} Let $A\in M_n(\p)$ satisfy ord $A\ge 2W+4m+2$. Then there exist $X\in M_n(\p)$ and $Y\in M_n(\p\{\tau\})$ (see the lines below (0.9.5)), $Y=\sum_{i=0}^{2m-1}Y_i\tau^i$, where $Y_i\in M_n(\p)$, such that:
\medskip
(1) $T_X(\goth B_2+Y)=(\goth B_2+Y)T_A$ (equality in $M_n(\p\{\tau\})$ );
\medskip
(2) $\goth B_2+Y\in GL_n(\p\{\tau\})$;
\medskip
(3) ord $Y_i\ge W_0$.
\medskip
(4) ord $X\ge \ord A-2m$.
\medskip
Moreover, there exists an algorithm of solution of 4.2 (1) (explained in the body of the proof) giving us a unique distinguished solution to 4.2 (1) (it satisfies to 4.2 (2), 4.2 (3) as well).
\medskip
{\bf Remark 4.3.} Conditions 4.2 (1), 4.2 (2) show that $M(A)$ is isomorphic to $M(X)$. Namely, in notations of (0.8) we denote $e_*$ for $M(A)$ (resp. $M(X)$) by $e_{*,A}$, resp. $e_{*,X}$. In these notations (0.8) is written as  $Te_{*,A}=T_Ae_{*,A}$, resp. $Te_{*,X}=T_Xe_{*,X}$, where multiplication by $T$ is in $\p[T,\tau]$ and multiplication by $T_A$, resp. $T_X$ is the action of $M_n(\p\{\tau\})$ on $\p\{\tau\}^{\oplus n}$. Let $\vf: M(X)\to M(A)$ be a $\p\{\tau\}$-linear map defined by the formula $$\vf(e_{*,X})=(\goth B_2+Y)e_{*,A}\eqno{(4.4)}$$ Condition 4.2 (1) is equivalent to the condition that $\vf$ is a $\p[T,\tau]$-homomorphism, and 4.2 (2) means that $\vf$ is a isomorphism.
\medskip
{\bf Proof of 4.2.} The equation (4.2 (1)) is
$$(\theta +X\tau+\tau^2)(\sum_{i=0}^mB_{1i}T_0^i+\sum_{i=0}^{2m-1}Y_i\tau^i)=(\sum_{i=0}^mB_{1i}T_0^i+\sum_{i=0}^{2m-1}Y_i\tau^i)(\theta +A\tau+\tau^2)\eqno{(4.5)}$$
We consider $A$, $B_{1i}$ as fixed, and we should solve (4.5) with respect to unknowns $X$, $Y_i$.
Slightly rewriting (3.2), we denote by $B_{2i}$ the coefficients in the equality $\sum_{i=0}^mB_{1i}T_0^i=\sum_{i=0}^{2m}B_{2i}\tau^i$, hence $B_{2i}=0$ for odd $i$, and for even $i$ we have $$B_{2i}=\sum_{j=i/2}^mk_{i/2,j}B_{1j}\eqno{(4.6)}$$

{\bf 4.7.} We have: ord$(B_{2i})\ge-(m-i/2)q^i$ ($i$ is even).
\medskip
 We consider the equality of coefficients of the terms at $\tau^{i+2}$ of (4.5), where $i=-1,\dots, 2m-1$. For the case $i=2m$ this equality is $B_{1m}^{(2)}=B_{1m}$ which is always satisfied because $B_{1m}\in M_n(\n F_{q^2})$; for $i=-2$ this equality is trivial. Using the fact that $T_0$ entering to the both $\theta +A\tau+\tau^2$, $\theta +X\tau+\tau^2$, commutes with $\sum_{i=0}^mB_{1i}T_0^i$ we see that this equality is the following:
$$Y_{i}^{(2)}+XY_{i+1}^{(1)}+XB_{2,i+1}^{(1)}+\theta Y_{i+2}=Y_{i}+Y_{i+1}A^{(i+1)}+ B_{2,i+1}A^{(i+1)}+ \theta^{q^{i+2}}Y_{i+2}\eqno{(4.8)}$$
(for $j\not\in \{0,1,\dots, 2m-1\}$ we have $Y_j=0$).

We shall use the following elementary
\medskip
{\bf Lemma 4.9.} Let $Z$ be an unknown matrix and $U$ a matrix parameter such that ord $U>0$. Then the equation $$Z=Z^{(2)}+U\eqno{(4.9.1)}$$ has a unique solution satisfying ord $Z>0$, it is given by the formula $$Z=U+U^{(2)}+U^{(4)}+U^{(6)}+...\eqno{\square \ \ (4.9.2)}$$

To solve (4.5) we consider first $X$ as a parameter, and we find $Y_i$ consecutively, $i=2m-1,\dots,0$, as functions of $X$. Equalities $(4.8)$ for $i=2m-1$, $2m-2,\ \dots,\ 0$ are equalities of the type (4.9.1) for the unknowns $Y_{i}$ if we treate $X$, $Y_{i+1}, \dots, Y_{2m-1}$ as parameters. As the last step, we substitute the values of $Y_0$, $Y_1$ (as functions of $X$) to the equality (4.8) for $i=-1$, and we solve it with respect to $X$. From the first sight, we can have a vicious circle: $Y_i$ are obtained as functions of $X$, and to find $X$ we should solve (4.8), $i=-1$; it contains $Y_i$.

Really, there is no vicious circle. Lemma 4.9 implies that formally
$$Y_i=D_i+D_i^{(2)}+D_i^{(4)}+D_i^{(6)}+...\eqno{(4.10)}$$
where $$D_i=XY_{i+1}^{(1)}+XB_{2,i+1}^{(1)}-\theta_{i+2,0} Y_{i+2}-Y_{i+1}A^{(i+1)}- B_{2,i+1}A^{(i+1)} \eqno{(4.11)}$$
Let $S_X(1)\subset M_n(\p)$ be the set of $X$ such that $D_{2m-1}$ of (4.11) satisfies ord $D_{2m-1}>0$ (we shall show later that for sufficiently small $A$ we have: $S_X(1)$ and subsequent $S_X(i)$, $S_{Y_j}(i)$ contain an open neighborhood of 0 in $M_n(\p)$). This means that $Y_{2m-1}$ --- the solution to (4.8), $i=2m-1$ --- can be considered as a function $f_{2m-1}: S_X(1) \to M_n(\p)$ defined by (4.10) for $i=2m-1$, i.e. if $Y_{2m-1}=f_{2m-1}(X)$ then (4.8), $i=2m-1$ holds for these $X$, $Y_{2m-1}$. Analogically, there exist sets $S_X(2)\subset S_X(1)$, $S_{Y_{2m-1}}(2)\subset M_n(\p)$ such that $S_{Y_{2m-1}}(2)\supset f_{2m-1}(S_X(2))$ and such that $Y_{2m-2}$ --- the solution to (4.8), $i=2m-2$ --- can be considered as a function $f_{2m-2}: S_X(2) \times S_{Y_{2m-1}}(2)\to M_n(\p)$ defined by (4.10) for $i=2m-2$, i.e. if $Y_{2m-2}=f_{2m-2}(X, f_{2m-1}(X))$ then (4.8), $i=2m-2$ holds for these $X$, $f_{2m-1}(X)$, $Y_{2m-2}$. We denote  $g_{2m-1}(X):=f_{2m-1}(X)$ and $g_{2m-2}(X):= f_{2m-2}(X, g_{2m-1}(X))$. Continuing the process, we get the sets $S_X(i)$, $S_{Y_j}(i)$, functions

$f_{2m-i}: S_X(i) \times S_{Y_{2m-(i-1)}}(i)\times S_{Y_{2m-(i-2)}}(i)\to M_n(\p)$ defined by (4.10) for $i=2m-i$ and functions
$g_{2m-i}$ defined by induction:

i.e. if $g_{2m-i}(X):=f_{2m-i}(X, g_{2m-(i-1)}(X), g_{2m-(i-2)}(X))$ then (4.8) for $2m-i$  holds for these $X$, $Y_j=g_{2m-j}(X)$, where $j=i, i-1,i-2$.

At the end of the process we get a function $f_0: S_X(2m)\times S_{Y_1}(2m)\times S_{Y_2}(2m) \to M_n(\p)$ such that  if $Y_0=g_0(X):=f_0(X, g_1(X), g_2(X))$ then (4.8) for $i=0$ holds. Obviously all entries of $g_i$ are power series on entries of $X$. Finally, we substitute $Y_0$, $Y_1$ that enter in $(4.8)$ for $i=-1$ by $g_0(X)$, $g_1(X)$ respectively. We get an equation $F(X)=0$ with one (matricial) unknown $X$ (and $A$, $B_{2i}$ parameters), where $F$ is a power series. If we choose a solution $X_0$ to the equation $F(X)=0$ and we let $Y_{0,i}:=g_i(X_0)$ then equations (4.8) for all $i$ are satisfied.
\medskip
Now we shall show that if $A$ is sufficiently small then

(1) All above $S_X(i)$, $S_{Y_j}(i)$ contain an open neighborhood of 0 in $M_n(\p)$, i.e. the condition ord $Z>0$ of 4.9 holds, and

(2) The equation $F(X)=0$ can be solved using Hensel lemma. This will give us a unique distinguished solution to (4.5).
\medskip
{\bf Lemma 4.12.} The formal expression for $Y_i$ obtained by the successive application of (4.9.2) is the following: $$Y_i=\sum\beta_{i,J,\Psi, l, \alpha} X^{(J)}\cdot B_{1l}^{(\alpha)} \cdot A^{(\Psi)}\eqno{(4.12.1)}$$ where the sum runs over all $J=(j_1,\dots,j_\mu)$, $\Psi=(\psi_1,\dots,\psi_\nu)$ satisfying $0\le j_1<j_2<...<j_\mu$ (i.e. $J$ is special in the terminology of the proof of Proposition 2), $\psi_1>\psi_2>...>\psi_\nu\ge i+1$ ($\mu$, $\nu$ can be 0, i.e. $J$, $\Psi$ can be $\emptyset$), $l$ satisfies $\frac{i+1}2\le l \le m$, $\alpha=0$ or 1, and $\beta_{i,J,\Psi, l, \alpha}\in \n F_q[\theta]$.
\medskip
{\bf Proof.} Induction by $i$ from $2m-1$ to $0$, using (4.11), (4.10). $\square$
\medskip
{\bf 4.13.} In principle, it can happen that applying (4.11), (4.10) for $i=2m-1, \dots, 0$, we can get for some given $J$, $\Psi$ more than one term of the form (4.12.1). Nevertheless, the same induction consideration shows that if this phenomenon occurs, we get only finitely many such terms with these $J$, $\Psi$. This means that we can consider $\beta_{i,J,\Psi, l, \alpha}$ as a finite sum of coefficents of the corresponding terms.
\medskip
After substitution of (4.12.1) for $Y_1$, $Y_0$ to (4.8), $i=-1$, we shall get an equality $F_{-1}(A,X)=0$ where
$$F_{-1}(A,X):=\sum\beta_{-1,J,\Psi, l, \alpha} X^{(J)}\cdot B_{1l}^{(\alpha)} \cdot A^{(\Psi)}\eqno{(4.14)}$$ is a series of the same type as (4.12.1). Like in (4.13), for any $J$, $\Psi$ there is no more than one such term. Let us evaluate coefficients $\beta_*$.
\medskip
{\bf Lemma 4.15.} Formulas (4.12.1), (4.14) do not contain terms having $\mu=\nu=0$.
\medskip
{\bf Proof.} Induction by $i$ from $2m-1$ to $-1$, using (4.11), (4.8) for $i=-1$. $\square$
\medskip
We shall call the terms of (4.12.1), (4.14) having $\mu=1$, $j_1=0$, $\nu=0$, resp. $\nu=1$, $\psi_1=0$, $\mu=0$ (i.e. terms of the form $\beta_*X B_{1l}^{(\alpha)}=X(\beta_*B_{1l}^{(\alpha)})$, resp. $(\beta_* B_{1l}^{(\alpha)})A$) the $X$- (resp. $A$)-principal terms, and we denote the right coefficient at $X$ (resp. the left one at $A$) of their sum by $P_{X,i}$, resp. $P_{A,i}$ ($i=2m-1,\dots,-1)$ (the $X$-, resp. $A$-principal part of $Y_i$ or of (4.14)).
\medskip
{\bf Lemma 4.16.} For $i$ even we have $P_{X,i}=0$. For $i$ odd, $\iota:=(i+1)/2$ we have
$$P_{X,i}=\overset{m}\to{ \underset{j=\iota}\to{\sum}}l_{\iota j}B_{1j}^{(1)}\eqno{(4.16.1)}$$ Particularly, for $i=-1$ we have $$P_{X,-1}=\overset{m}\to{ \underset{j=0}\to{\sum}}\theta^jB_{1j}^{(1)}\eqno{(4.16.2)}$$
Further, $P_{A,i}=0$ for $i\ge0$, and $P_{A,-1}=\overset{m}\to{ \underset{j=0}\to{\sum}}\theta^jB_{1j}$.
\medskip
{\bf Proof.} Induction by $i$ from $2m-1$ to $-1$. The only terms of (4.10), (4.11) that contribute to the $X$-principal part of $Y_i$ is the term $D_i$ of (4.10) and the terms $XB_{2,i+1}^{(1)}$, $-\theta_{i+2,0} Y_{i+2}$ of (4.11), because Lemma 4.15 shows that the term $XY^{(1)}_{i+1}$ of (4.11) cannot contribute to the principal part. Since $ B_{2i}=0$ for odd $i$, we get immediately that $P_{X,i}=0$ for even $i$. This gives us the recurrence relation:
$$P_{X,i}=B_{2,i+1}^{(1)}-\theta_{i+2,0}P_{X,i+2}\eqno{(4.16.3)}$$ with the initial condition $P_{X,2m+1}=0$. (4.6), (4.16.3) and Lemma 3.6 ($i$ of Lemma 3.6 is $\iota$ of the present lemma) give us immediately (4.16.1).
\medskip
(4.11) shows that $P_{A,i}=0$ for $i\ge0$. (4.8) for $i=-1$ shows that the terms that can contribute to the $A$-principal part of (4.14) are $Y_{0}A$, $B_{20}A$. Lemma 4.15 shows that $Y_{0}A$ do not contribute to the $A$-principal part of (4.14), hence $P_{A,-1}=B_{20}=\overset{m}\to{ \underset{j=0}\to{\sum}}\theta^jB_{1j}$. $\square$
\medskip
Let us evaluate ord $\be_*$ from (4.12.1), (4.14). For $J$, $\Psi$ of (4.12.1), (4.14) we denote $\goth m=\goth m(J,\Psi):= \max (j_*, \psi_*)$ (if $J$, resp. $\Psi=\emptyset$, then we let max $j_*$, resp. max $\psi_*=0$).
\medskip
{\bf Lemma 4.17.} For all terms of (4.12.1), (4.14) we have
$$\ord (\beta_*) + W q^{\goth m(J,\Psi)}\ge W_0\eqno{(4.17.0)}$$

{\bf Proof.}  Induction by $i$ from $2m-1$ to $-1$ (here and below the case $i=-1$ means the formula (4.14)). More exactly, we prove by induction that for all terms of $Y_i$ we have
$$\ord (\beta_*) + W q^{\goth m(J,\Psi)}\ge (i+2)W_0\eqno{(4.17.1)}$$
If (4.17.1) holds for all terms that enter in $D_i$ of (4.10) then it holds for terms $D_i^{(2j)}$ as well, hence it is sufficient to prove (4.17.1) only for $D_i$. Further,
\medskip
{\bf (4.17.2)} All number coefficients that enter in (4.11), namely $B_{2,i+1}$, $B_{2,i+1}^{(1)}$ and $\theta_{i+2,0}$, have ord $\ge -W_0$ (see (3.3)).
\medskip
For $i=2m-1$ the formula (4.11) becomes $D_{2m-1}=XB_{2,2m}^{(1)}-\goth B_{2,2m}A^{(2m)}$ and (4.17.1) obviously holds. The induction step from $i$ to $i-1$ also follows immediately from (4.17.2). $\square$
\medskip
{\bf Lemma 4.18.} If $\ord X>W$, $\ord A>W$ then the series $Y_i=Y_i(X,A)$, $F_{-1}=F_{-1}(X,A)$ converge.
\medskip
{\bf Proof.}  Let $\ve>0$ satisfies $\ord X>W+\ve$, $\ord A>W+\ve$. Ord of a term of (4.12.1), (4.14) is $> \ord(\beta_*) + (W+\ve) q^\goth m> \ve q^\goth m+W_0 $. Since for any given $\goth m_0$ there are only finitely many terms in (4.12.1), (4.14) having $\goth m\le\goth m_0$ we get that ords of terms of (4.12.1), (4.14) tend to $+\infty$. $\square$
\medskip
Recall that $\goth B_1=\overset{m}\to{ \underset{j=0}\to{\sum}}B_{1j}\theta^j$; we denote $\bar \goth B_1 =\overset{m}\to{ \underset{j=0}\to{\sum}}B_{1j}^{(1)}\theta^j$, i.e. bar means the $\n F_{q^2}/\n F_q$-conjugation.
\medskip
{\bf Lemma 4.19.} The only terms of $F_{-1}(A,X)$ having $\goth m=0$ are $X$- and $A$-principal terms, i.e. $X\bar \goth B_1$ and $-\goth B_1A$.
\medskip
{\bf Proof.} We must prove that $F_{-1}(A,X)$ does not contain terms of the form $\beta_*XB_{1l}^{(\alpha)}A$. These terms can appear only from the term $Y_{i+1}A^{(i+1)}$ of (4.11) for $i=-1$. But $P_{X,0}=0$. $\square$
\medskip
{\bf Proposition 4.20.} Let $W_1:=2W+2m+2$. Let $A$ be such that $\ord A\ge W_1$. Then there exists $X$ having $\ord X\ge \ord A-2m$ which is a root to the equation $F_{-1}(A,X)=0$.
\medskip
{\bf Proof.}  We shall show that the conditions of Lemma 2.29 hold for this case. We multiply $F_{-1}(A,X)=0$ by $\bar \goth B_1^{-1}$ from the right, we get
$$X=\goth B_1A\bar \goth B_1^{-1} + \hbox{ (terms having $\goth m\ge 1$) }\cdot\bar \goth B_1^{-1}\eqno{(4.20.1)}$$
Since $\ord \goth B_1, \ord \bar \goth B_1^{-1}\ge-m$ we have that $\ord \goth B_1A\bar \goth B_1^{-1}\ge 2W+2$. We change the scale $X=\de X'$ where ord $\de=2W+1$. The equation $F_{-1}(A,X')=0$ is of the form (2.28). Let us evaluate ord $u$ and ord $C(J,i)$.
\medskip
Case 1: $j_\mu\ge \psi_1$ or $\Psi=\emptyset$. In this case we have $\goth m(J,\Psi)=j_\mu\ge m(J)-1$, $j_\mu\ge1$,
$$\ord C(J,i)\ge \min_{\Psi,l,\al}(\ord \be_{-1,J,\Psi,l,\al})+(2W+1)q^{j_\mu}-(2W+1)-m\ge j_\mu+1\ge m(J)$$
(for the first inequality, we take into consideration only the highest power $X^{(j_\mu)}=\de^{q^{j_\mu}}{X'}^{(j_\mu)}$, and we neglect a  possible factor $A^{(\Psi)}$; for the second inequality, we use 4.17.0). Hence, for this case 2.29.2 is satisfied for $\ga=1$.
\medskip
Case 2: $j_\mu< \psi_1$ and $J\ne\emptyset$. In this case we have $\goth m(J,\Psi)=\psi_1\ge m(J)$. A term $C(J,i)$ corresponding to $(-1, J,\Psi,l,\al)$ in (4.14) for some $l,\al$ satisfies
$$\ord C(J,i)\ge \ord \be_{-1,J,\Psi,l,\al}+(2W+2m+2)q^{\psi_1}-(2W+1)-m\ge \psi_1\ge m(J)\eqno{(4.20.2)}$$
(for the first inequality, we take into consideration only the highest power $A^{(\psi_1)}$ whose ord $\ge(2W+2m+2)q^{\psi_1}$, and we neglect a  possible factor $\de^{(J)}$; for the second inequality, we use 4.17.0). Hence, for this case 2.29.2 is satisfied as well for $\ga=1$.

Further, since for given $J$, $\Psi$ there is no more than one term of (4.14) having these $J$, $\Psi$, for a given $J$, $\psi_1$ there are only finitely many terms of (4.14) having these $J$, $\psi_1$. (4.20.2) shows that 2.29.3 is satisfied. Finally,  $u=\de^{-1}\goth B_1A\bar \goth B_1^{-1}+$ the sum of other terms having $J=\emptyset$. The first two inequalities of (4.20.2) hold for $J=\emptyset$, hence the series for $u$ converges and ord $u\ge\ord A-(2W+1)-2m\ge1$, hence 2.29.1 is satisfied. Since $\ord X'=\ord u$ we get that $\ord X\ge \ord A-2m$. $\square$
\medskip
This proves 4.2(1). Condition 4.2(3) is obviously satisfied, because of 4.12.1, 4.17.1 and inequalities ord $A\ge W$, ord $X\ge W$. Condition 4.2(4) also follows from the above. Let us prove 4.2(2).
\medskip
{\bf Proposition 4.21.} If ord $A\ge W_1+2m$ then $\goth B_2+Y$ is invertible, i.e. 4.2 (2) holds.
\medskip
{\bf Proof.} We call $A$, $\goth B_1$ as the input data of Proposition 4.2, and $X$, $Y$ as the output data. We have ord $X\ge W_1$, hence we can apply 4.2(1) for the input data $X$, $\goth B_1^{-1}$. We denote the corresponding output data as $X'$, $Y'$.
We shall get the equality 4.2.1 for this situation
$$(\theta +X'\tau+\tau^2)(\goth B^{-1}_2+Y')=(\goth B^{-1}_2+Y')(\theta +X\tau+\tau^2)$$
Multiplying it by $\goth B_2+Y$ from the right and using (4.2.1) we get
$$(\theta +X'\tau+\tau^2)(\goth B^{-1}_2+Y')(\goth B_2+Y)=(\goth B^{-1}_2+Y')(\goth B_2+Y)(\theta +A\tau+\tau^2)\eqno{(4.21.1)}$$
We have $(\goth B^{-1}_2+Y')(\goth B_2+Y)=1+\De$ where $\De=\goth B^{-1}_2Y+Y'\goth B_2+Y'Y$. (4.7) and the condition ord $Y$, ord $Y'>W_0$ imply ord $\De>0$. Let $\De_i\tau^i$ be the highest non-0 term of $\De$ as a polynomial in $\tau$. The equality of $\tau^{i+2}$-terms of both sides of (4.21.1) is $\De_i^{(2)}=\De_i$, which contradicts to the condition ord $\De>0$. This means that $(\goth B^{-1}_2+Y')(\goth B_2+Y)=1$. It is well-known that in this case $(\goth B_2+Y)(\goth B^{-1}_2+Y')=1$ as well. Really, the ring $\p\{\tau\}$ satisfies both (left and right) Ore condition ([G], Section 1). Since it is without zero-divisors, this implies that $\p\{\tau\}$ is a subring of a skew field $K$, see, for example, [C], Theorem 1.2.2, p. 7-8. Finally, there is a theorem that if an element of $M_n(K)$ has a right (or left) inverse then it is invertible, see, for example, [D], Section 19, Theorem 3, p. 131 (it is clear that $(\goth B_2+Y)^{-1}$ which a priori belongs to $M_n(K)$ is really $\goth B^{-1}_2+Y'$). $\square$
\medskip
Therefore, Proposition 4.2 is proved. $\square$
\medskip
To pass from t-motives to lattices, we consider $V_1$, $V_2$, $L_1$, $L_2$ from (1.6). Let $f_1,\dots, f_n$, $e_1,\dots, e_n$ be some bases of $V_1$, $V_2$ respectively. We denote the vector columns $(f_1,\dots, f_n)^t$, $(e_1,\dots, e_n)^t$ by $f_*$, $e_*$ respectively. Let $\vf: V_1\to V_2$ from (1.6) (particularly, $\vf(L_1)=L_2$ ) be defined by the formula $$\vf(f_*)=F^t\cdot e_*\eqno{(4.22)}$$  where $F\in M_n(\p)$. We consider the following situation (notations of (1.6)): $\goth g_1$ (resp. $\goth h_1$) is close to $f_*$, resp. $e_*$, and $\goth g_2$ (resp. $\goth h_2$) is close to $\omega f_*$, resp. $\omega e_*$ (see statement of Lemma 4.23 for the exact estimates).

For $x=(c_1,\dots,c_n)\in \p^n$, where $c_i\in \p$, we let (like for the square matrices) ord$(x)=\underset{i}\to{\hbox{ min }}\hbox{ ord }(c_i)$. For a set of vectors $x_*=(x_1,\dots,x_n)^t$, where $x_i\in \p^n$, we let ord$(x)=\underset{i}\to{\hbox{ min }}\hbox{ ord }(x_i)$. For two bases (for example, $f_*$ and $\goth g_1$) the distance between them is given by ord $f_*-\goth g_1$, where $f_*-\goth g_1:=(f_1-g_1,\dots,f_n-g_n)$. Analogically, we define the distance between 2 matrices in $M_n(\p)$.

Finally, let $F$ of (4.22) satisfy ord $F\ge -m$ and let there exists $D\in G_1$ such that ord $(F-D)>0$. Let $Z$ be from (1.6.1).
\medskip
{\bf Lemma 4.23.} In the above notations let we have ord $(f_*-\goth g_1)>m$, ord $(\omega f_*-\goth g_2)>m$,  ord $(e_*-\goth h_1)>m$, ord $(\omega e_*-\goth h_2)>m$. Then $Z\in G_0$, and $\beta(Z)=D$.
\medskip
{\bf Proof.} Let $\Delta_1$, $\Delta_2$, $\Delta_3$, $\Delta_4$, $Y_0$ be matrices defined by the conditions $\goth g_1=f_*+\Delta_1 f_*$, $$\goth h_1=e_*+\Delta_2 e_*, \ \ \goth h_2=\omega e_*+\Delta_3 e_*\eqno{(4.23.0)}$$ $\goth g_2=\omega f_*+\Delta_4 f_*$, $F=D+Y_0$ respectively. Further, the condition $D\in G_1$ implies $D=D_1+\omega D_2$, where $D_1$, $D_2\in M_n(\n F_q[\theta])$.  We have $\vf(\goth g_1)= F^te_*+\Delta_1F^te_*=(I_n+\Delta_1)F^t(I_n+\Delta_2)^{-1}\goth h_1$, hence
\medskip
$\vf(\goth g_1)=(I_n+\Delta_1)(D_1^t+Y_0^t)(I_n+\Delta_2)^{-1} \goth h_1+ (I_n+\Delta_1)D_2^t (I_n+\Delta_2)^{-1} \omega\goth h_1$.
\medskip
Since $\goth h_1=(I_n+\Delta_2)(\omega I_n+\Delta_3)^{-1}\goth h_2$, we get $$\vf(\goth g_1)=(I_n+\Delta_1)(D_1^t+Y_0^t)(I_n+\Delta_2') \goth h_1+(I_n+\Delta_1)D_2^t(I_n+\Delta_3')\goth h_2\eqno{(4.23.1)}$$ where $\Delta_2', \ \Delta_3'$ are defined by the conditions $I_n+\Delta_2'=(I_n+\Delta_2)^{-1}$, $I_n+\Delta_3'=\omega(\omega I_n+\Delta_3)^{-1}$. Since $\vf(L_1)=L_2$, we get that all elements of the vector column (4.23.1) (they are elements of $V_2$ ) belong to $L_2$. From another side, all elements of the vector column $D_1^t \goth h_1+D_2^t \goth h_2$ also belong to $L_2$. Obviously we have for their difference:
$$(I_n+\Delta_1)(D_1^t+Y_0^t)(I_n+\Delta_2') \goth h_1+(I_n+\Delta_1)D_2^t(I_n+\Delta_3')\goth h_2- (D_1^t \goth h_1+D_2^t \goth h_2)\in (L_2)^n\eqno{(4.23.2)}$$
and
$$(I_n+\Delta_1)(D_1^t+Y_0^t)(I_n+\Delta_2') \goth h_1+(I_n+\Delta_1)D_2^t(I_n+\Delta_3')\goth h_2- (D_1^t \goth h_1+D_2^t \goth h_2)=$$ $$=\Delta_5\goth h_1+\Delta_6\goth h_2, \ \  \hbox{ where $\Delta_5, \ \Delta_6\in M_n(\p)$ and ord } \Delta_5>0, \hbox{ ord } \Delta_6>0.\eqno{(4.23.3)}$$
It is easy to see that (4.23.2), (4.23.3) imply $\Delta_5=\Delta_6=0$. Really, considering the $i$-th line of $\Delta_5, \ \Delta_6$ and respectively the $i$-th component of $\Delta_5\goth h_1+\Delta_6\goth h_2$, we see that it is sufficient to prove an obvious
\medskip
{\bf Sublemma 4.23.4.} Let we have an equality $$\sum_{i=1}^{2n} \delta_{7i}h_i=\sum_{i=1}^{2n} \alpha_{i}h_i\eqno{(4.23.4.1)}$$ where $\delta_{7i}\in \p$, ord $\delta_{7i}>0$ and $\alpha_i\in \n F_q[\theta]$. Then $\alpha_{i}=0$.
\medskip
{\bf Proof.} Let, conversely, $\alpha_{i}\ne0$. We denote $k=$ max deg $\alpha_i$ (here $\alpha_i$ are treated as polynomials in $\theta$). We multiply the equality (4.23.4.1) by $\theta^{-k}$: $$\sum_{i=1}^{2n} \theta^{-k}\delta_{7i}h_i=\sum_{i=1}^{2n} \theta^{-k}\alpha_{i}h_i\eqno{(4.23.4.2)}$$The coefficients of the left hand side $\theta^{-k}\delta_{7i}$ have ord $> 0$, the coefficients of the right hand side $\theta^{-k}\alpha_{i}$ have ord $\ge0$ and at least one of these coefficients has ord = 0. Conditions (4.23.0) imply that $\forall i$ coordinates of $h_i$ in the basis $e_*$ have ord $\ge0$, i.e. they belong to $O_{\p}$, hence we can consider the reduction of (4.23.4.2) (treated as a line of coordinates) modulo $\goth m$ --- the maximal ideal of $O_{\p}$ ( = the set of elements having ord $> 0$). Again (4.23.0) imply that the reduction of $h_i$ is equal to $e_i$, resp. $\omega e_i$ for $i\le n$, resp. $i>n$. The reduction of the left hand side of (4.23.4.2) is 0, while the reduction of  the right hand side of (4.23.4.2) is $\sum_{i=1}^{n}\widetilde{\theta^{-k}\alpha_{i}} \ e_i+\sum_{i=1}^{n}\widetilde{\theta^{-k}\alpha_{n+i}} \ \omega e_i$. Since all coefficients $\widetilde{\theta^{-k}\alpha_{i}}, \ \widetilde{\theta^{-k} \alpha_{n+i}} \in \n F_q$ and not all of them are 0, we get a contradiction.  $\square$
\medskip
So, we get that $\vf(\goth g_1)=D_1^t \goth h_1+D_2^t \goth h_2$. The proof that $\vf(\goth g_2)=kD_2^t \goth h_1+D_1^t \goth h_2$ is completely analogous. $\square$
\medskip
{\bf 4.24. End of the proof.} The main theorem follows easily from Proposition 4.2 and Lemma 4.23. We define $U_m$ as the set of $A$ such that ord $A\ge 2W+4m+2$, as in (4.2). Let $A_1$, $A_2$, $\gamma$ be from (0.9). We consider the Proposition 4.2 for the case $A=A_1$, $\goth B_1=\beta(\gamma)$, $\goth B_2=\alpha\circ\beta(\gamma)$. Let $\vf:M(X)\to M(A)$ be the isomorphism of Remark 4.2.1. It  remains to prove that $X=A_2$. Since obviously both $A_2, \ X\in \goth W_1$ (see beginning of Section 2) and $\Cal S$ is an isomorphism on $ \goth W_1$, we have: $X=A_2$ is equivalent to the condition $\Cal S(X)=\Cal S(A_2)$, and, since $\Cal S(A_2)=\ga(\Cal S(A_1))$, it is sufficient to prove that $\Cal S(X)=\ga(\Cal S(A_1))$.

Let us consider T-modules $E(A)$, $E(X)$ corresponding to $M(A)$, $M(X)$ respectively (see [G], Definition 5.4.5 and below). The functor $E$ is contravariant, hence we have a map $E(\vf): E(A) \to E(X)$. Both $E(A)$, $E(X)$ are isomorphic to $\p^n$. We identify elements of $E(A)$, $E(X)$ with $n\times1$-matrix columns of elements of $\p$. (4.2.1.1) implies that in this matrix form $E(\vf)$ is given by the formula $E(\vf)(C)=\sum_{i=0}^{2m} ( B_{2i}+Y_i) C^{(i)}$ (because $\goth B_2+Y$ of (4.2.1.1) $=\sum_{i=0}^{2m} (B_{2i}+Y_i)\tau^i$), where $C\in E(A)=\p^n$.

Let us consider the induced Lie morphism $E(\vf)_*: \Lie(E(A))\to \Lie(E(X))$ (see [G], below Definition 5.9.4, p. 161, lines 2 - 3 from the bottom). We identify $\Lie(E(A))$ and $\Lie(E(X))$ with $\p^n$. It follows from [G], a formula above Remark 5.9.8, that $E(\vf)_*(C)=(B_{20}+Y_0)C=(\goth B_1+Y_0)C$,  where $C\in \Lie(E(A))=\p^n$.

Now we apply Lemma 4.23 for the case $V_1=\Lie(E(A))$, $V_2=\Lie(E(X))$, $\vf=E(\vf)_*$, $L_1$, $L_2$ are the lattices of $M(A)$, $M(X)$ respectively, $f_*=e_*=\goth e$ (Section 2, below 2.7), $\goth g_i=\goth g_i(A)$, $\goth h_i=\goth g_i(X)$, where $i=1,2$, $\goth g_i(*)$ are from Section 2, below 2.9. We see that $F$ from (4.22) is $\goth B_1+Y_0$. This means that $D$ from 4.23 is $\goth B_1$. Inequalities on ord $Y$ (Section 4), ord $\Delta$, ord $\Delta'$ (Section 2) show that all conditions of Lemma 4.23 are satisfied, hence $U$ from (1.6.1) is equal to $\be^{-1}(\goth B_1)=\ga$, hence $\Cal S(X)=\ga(\Cal S(A_1))$. This proves the theorem 0.9.
\medskip
{\bf  Remark 4.25.} Let us show the concordance of formulas for the principal term of $X$ as a function of $A$, $\ga$. From one side, we have a formula (4.20.1). From another side, (2.26) shows that for some $\la\ne0$ we have $\Cal S(A)=\omega+\lambda A^t+$ higher terms, $S(X)=\omega+\lambda X^t+$ higher terms. These formulas and Lemma 4.1 show that $A^t=\bar \goth B_1^tX^t\goth B_1^{t-1}$ + higher terms which is equivalent to (4.20.1).
\medskip
{\bf References}

\nopagebreak
\medskip
[A] Anderson Greg W. $t$-motives. Duke Math. J. Volume 53, Number 2 (1986), 457-502.
\medskip
[C] Cohn, P. M. Skew fields. Theory of general division rings. Encyclopedia of Mathematics and
its Applications, 57. Cambridge
University Press, Cambridge, 1995. xvi+500 pp.
\medskip
[D] Draxl, P. K. Skew fields.
London Mathematical Society Lecture Note Series, 81.
Cambridge University Press, Cambridge, 1983. ix+182 pp.
\medskip
[Dr]  Drinfeld, V. G. Elliptic modules. Math. USSR-Sb. no. 4, 561 - 592 (1976).
\medskip
[G] Goss, David Basic structures of function
field arithmetic.
Springer-Verlag, Berlin, 1996. xiv+422 pp.
\medskip
[H] Urs Hartl, Uniformizing the Stacks of Abelian Sheaves.

http://arxiv.org/abs/math.NT/0409341
\medskip
[L] Logachev, Duality of Anderson t-motives. arxiv.org/pdf/0711.1928.pdf
\medskip
[L1] Logachev, Anderson t-motives are analogs of abelian varieties with multiplication by imaginary quadratic fields. arxiv.org/pdf/0907.4712.pdf
\medskip
[L2] Logachev, Fiber of the lattice map for t-motives. In preparation.

\medskip
E-mails: shuragri{\@}gmail.com, logachev94{\@}gmail.com

\enddocument